\documentclass[a4paper,10pt]{amsart}

\newtheorem{theo}{Theorem}[section]

\newtheorem{prop}[theo]{Proposition}
\newtheorem{lem}[theo]{Lemma}
\newtheorem{cor}[theo]{Corollary}

\newtheorem{rema}[theo]{Remark}
\newtheorem{remas}[theo]{Remarks}

\def \kbar {{\bar k}}

\def \dem {\paragraph{ \em Proof. }}

\def \Romannumeral #1 {\expandafter\uppercase\expandafter {\romannumeral #1} }
\def \br {{\rm{Br \,}}}

\def \P {{\bf P}}

\def \pic {{\rm {Pic\,}}}
\def \Div {{\rm{Div\,}}}
\def \gal {{\rm{Gal\,}}}

\def \calo {{\mathcal O}}
\def \T {{\mathcal T}}
\def \X {{\mathcal X}}
\def \D {{\mathcal D}}

\def \calf {{\mathcal F}}
\def \calg {{\mathcal G}}
\def \calt {{\mathcal T}}
\def \spec {{\rm{Spec\,}}}

\def \dim {{\rm{dim\,}}}

\def \Hom {{\rm {Hom}}}
\def \Ext {{\rm {Ext}}}

\def\ov{\overline}

\def \Z {{\bf Z}}
\def \Q {{\bf Q}}

\def \RR {{\bf R}}

\def \im {{\rm {Im}}}
\def \ker {{\rm {Ker}}}
\def \G {{\bf G}_m}
\def \A {{\bf A}}

\def\smallsquare{\vbox{\hrule\hbox{\vrule height 1 ex\kern 1 ex\vrule}\hrule}}
\def\enddem{\hfill \smallsquare\vskip 3mm}

\def \abstract{\paragraph{Abstract. }}
\usepackage{amscd}
\usepackage{amssymb}
\usepackage{amsmath}

\everymath{\displaystyle}

\title{Local-global questions for tori over $p$-adic function fields}
\author{David Harari and Tam\'as Szamuely}
\address{Universit\'e de Paris-Sud Math\'ematique, B\^atiment 425, 91405 Orsay, France}
\email{David.Harari@math.u-psud.fr}
\address{Alfr\'ed R\'enyi Institute of Mathematics, Hungarian Academy of Sciences, Re\'altanoda utca 13--15, H-1053 Budapest, Hungary and Central European University, N\'ador utca 9, H-1051 Budapest, Hungary}
\email{szamuely.tamas@renyi.mta.hu}

\date{\today}

\def \id {{\rm id}}

\DeclareFontFamily{U}{wncy}{}
\DeclareFontShape{U}{wncy}{m}{n}{%
   <5>wncyr5%
   <6>wncyr6%
   <7>wncyr7%
   <8>wncyr8%
   <9>wncyr9%
   <10>wncyr10%
   <11>wncyr10%
   <12>wncyr6%
   <14>wncyr7%
   <17>wncyr8%
   <20>wncyr10%
   <25>wncyr10}{}
\DeclareMathAlphabet{\cyrille}{U}{wncy}{m}{n}
\def\Sha{\cyrille X}
\def \R{{\bf R}}

\def \Xbar{{\overline X}}

\usepackage{graphicx}

\makeatletter
\DeclareRobustCommand\widecheck[1]{{\mathpalette\@widecheck{#1}}}
\def\@widecheck#1#2{%
    \setbox\z@\hbox{\m@th$#1#2$}%
    \setbox\tw@\hbox{\m@th$#1%
       \widehat{%
          \vrule\@width\z@\@height\ht\z@
          \vrule\@height\z@\@width\wd\z@}$}%
    \dp\tw@-\ht\z@
    \@tempdima\ht\z@ \advance\@tempdima2\ht\tw@ \divide\@tempdima\thr@@
    \setbox\tw@\hbox{%
       \raise\@tempdima\hbox{\scalebox{1}[-1]{\lower\@tempdima\box
\tw@}}}%
    {\ooalign{\box\tw@ \cr \box\z@}}}
\makeatother

\begin{document}

\maketitle

\noindent{\small {\bf Abstract.} We study local-global questions for Galois cohomology over the function field of a curve defined over a $p$-adic field (a field of cohomological dimension 3). We define Tate-Shafarevich groups of a commutative group scheme via cohomology classes locally trivial at each completion of the base field coming from a closed point of the curve. In the case of a torus we establish a perfect duality between the first Tate-Shafarevich group of the torus and the second Tate-Shafarevich group of the dual torus. Building upon the duality theorem, we show that the failure of the local-global principle for rational points on principal homogeneous spaces under tori is controlled by a certain subquotient of a third \'etale cohomology group. We also prove a generalization to principal homogeneous spaces of certain reductive group schemes in the case when the base curve has good reduction.}

\setcounter{section}{-1}

\section{Introduction}

In recent years there has been considerable interest in local-global principles for group schemes defined over fields of cohomological dimension strictly greater than 2. In particular, Harbater, Hartmann and Krashen, in a series of papers of which perhaps \cite{hhk} is the most relevant to our present context, have used patching techniques to introduce and study an analogue of the Tate--Shafarevich group for certain linear algebraic groups  over the function field of a curve defined over a complete discretely valued field. In \cite{ctps} and the recent preprint \cite{ctps2}, Colliot-Th\'el\`ene, Parimala and Suresh have obtained local-global results over similar fields whose formulation shows a closer analogy with the classical theory over number fields.

These developments motivate a systematic study of local-global questions for Galois cohomology over the function field $K$ of a curve defined over a finite extension of $\Q_p$. In the present paper, except for the last section, we limit ourselves to tori over $K$ and prove two kinds of results: duality theorems for Galois cohomology, and analogues of the by now classical results of Sansuc \cite{sansuc} regarding the Hasse principle for torsors under tori. In the last section we also prove a generalization to torsors under certain reductive groups. Concerning duality, a direct precursor of our theorems is the Artin--Verdier style statement of Scheiderer and van Hamel \cite{schvh}.

Let us now review the main results of the paper. Given a torus $T$ over a field $K$ as above, set
$$
\Sha^i(T):=\ker[H^i(K,T)\to\prod_{v\in X^{(1)}}H^i(K_v,T)]
$$
where the product is over closed points of the smooth proper curve $X$ whose function field is $K$, and $K_v$ is the completion of $K$ with respect to the discrete valuation coming from the point $v$. Note that this definition is different from that of \cite{ctps} and \cite{ctps2} where all discrete rank 1 valuations of $K$ are considered, and also from that of \cite{hhk} where the authors complete the local rings of closed points of the special fibre of an integral model of $K$.

\begin{theo}\label{vardualsha}{\rm (= Theorem \ref{dualsha})} With notations as above, there is a perfect pairing of finite
groups
$$\Sha^1(T) \times \Sha^2(T') \to \Q/\Z.$$
\end{theo}
Here $T'$ is the dual torus of $T$, i.e. the torus whose character group is the cocharacter group of $T$. As we shall explain in Section \ref{one}, the duality pairing comes from a cup-product associated with a natural map $T\otimes T'\to \Z(2)[2]$ in the derived category, where $\Z(2)$ is the \'etale motivic complex of weight 2.

As an application of Theorem \ref{vardualsha}, we study cohomological obstructions to the Hasse principle over $K$. Given a smooth $K$-variety $Y$, its set $Y(\A_K)$ of adelic points consists of families $(P_v)$ with $P_v\in Y(K_v)$ for $v\in X^{(1)}$ such that $P_v$ is integral with respect to the ring of integers $\calo_v\subset K_v$  for almost all $v$. The Hasse principle is said to hold for $Y$ if $Y(\A_K)\neq\emptyset$ implies $Y(K)\neq\emptyset$. If $T$ is a $K$-torus with $\Sha^1(T)\neq 0$ (such examples exist; see Remark \ref{exsha} below), a nontrivial class in $\Sha^1(T)$ corresponds to a $K$-torsor that violates the Hasse principle.

In the classical case, obstructions to the Hasse principle are associated with subgroups of the Brauer group, following a fundamental insight of Manin. Here we are dealing with a base field of cohomological dimension 3, and therefore a third \'etale cohomology group plays the role of the Brauer group, namely $H^3(Y, \Q/\Z(2))$ (where $\Q/\Z(2)$ denotes the direct limit of the torsion \'etale sheaves $\mu_n^{\otimes 2}$ for all $n$). As explained in detail in Section \ref{five}, evaluating classes at adelic points of a $K$-variety $Y$ induces a pairing
$$
Y(\A_K)\times H^3(Y, \Q/\Z(2))\to\Q/\Z
$$
which annihilates the diagonal image of $Y(K)$ in $Y(\A_K)$ on the left, as well as
the image of the natural map $H^3(K,\Q/\Z(2))\to H^3(Y,\Q/\Z(2))$
on the right. The induced pairing on the group
\begin{align*}
H^3_{\rm lc}(Y,\Q/\Z(2)) &:= \ker[H^3(Y,\Q/\Z(2))/\im(H^3(K,\Q/\Z(2)))\to \\ &\to\prod_vH^3(Y\times_K{K_v},\Q/\Z(2))/\im(H^3(K_v,\Q/\Z(2)))]\end{align*}
is constant, whence a homomorphism
$$
\rho_Y:\,  H^3_{\rm lc}(Y,\Q/\Z(2))\to \Q/\Z.
$$
We can now state:

\begin{theo}\label{varmainph} {\rm (= Theorem \ref{mainph})}
Assume that $Y$ is a torsor under a $K$-torus $T$ such that $Y(\A_K)\neq\emptyset$ and $\rho_Y$ is the zero map. Then $Y(K)\neq\emptyset$.
\end{theo}

 In fact, the proof will show that it is enough to work with the image of a certain map $\Sha^2(T')\to H^3_{\rm lc}(Y,\Q/\Z(2))$ that relates the map $\rho$ to the duality pairing of Theorem \ref{vardualsha}. This phenomenon is analogous to what happens in the proofs of the corresponding facts over number fields, but there are considerable additional technical difficulties arising from higher cohomological dimension.

Finally, in the last section we leave tori and consider general reductive groups. Combining Theorem \ref{varmainph} with known facts concerning the Rost invariant for simply connected linear algebraic groups, we shall prove:

\begin{theo} {\rm (= Theorem \ref{phlag})}
Let $E$ be a torsor under a reductive linear algebraic group $G$ over $K$. Assume that the smooth proper  curve $X$ whose fraction field is $K$ has good reduction and that the simply connected cover $G^{\rm sc}$ of the derived subgroup of $G$  is quasi-split and has no $E_8$ factor. If $E(\A_K)\neq\emptyset$ and $\rho_E$ is the zero map, then $E(K)\neq\emptyset$.
\end{theo}

As in classical proofs, we exploit the technique of $z$-extensions first used by Deligne and Langlands to reduce to the cases of simply connected groups and tori. However, a remarkable difference with the situation over number fields is that here over the completions $K_v$ there may exist nontrivial torsors under $G^{\rm sc}$.

In the joint paper \cite{pcwa} with Claus Scheiderer that complements the present one, we construct a 9-term Poitou--Tate type exact sequence for tori over $K$ (and also a 12-term sequence for finite modules). Like in the number field case, part of the sequence can then be used to analyze the defect of weak approximation for a torus. We also show that the defect of weak approximation is controlled by a certain subgroup of the third unramified cohomology group of the torus, in analogy with Theorem \ref{varmainph}.

 We would like to point out that Theorem \ref{vardualsha} was announced by Claus Scheiderer in lectures given in 2002 but he did not publish his proof. Also, his intended approach was slightly different (see also Remark \ref{avrema} below).
We are very grateful to Joseph Ayoub, Michel Brion, Cyril Demarche, Mathieu Florence, Philippe Gille and above all Jean-Louis Colliot-Th\'el\`ene for instructive discussions and their help with several proofs. Part of this work has been done while the authors were visiting each other's home institutions whose hospitality was greatly appreciated. The first author would also like to thank the Centre Bernoulli (Lausanne) for hospitality and the second author the Hungarian Science Foundation OTKA for partial support under grant No. NK81203.

\bigskip

\noindent{\bf Notation and conventions.} The following conventions will be adopted throughout the article.\smallskip

\noindent{\em Abelian groups.\/}  Given an abelian group $A$, we shall denote by $\overline A$
the quotient of $A$ by its maximal divisible subgroup $\Div A$. The $n$-torsion subgroup and (for $\ell$ prime) the $\ell$-primary torsion subgroup of  $A$ will be denoted by $_n A$ and $A \{ \ell\}$, respectively.
If $A$ is torsion, it is said to be {\it of cofinite type} if for each
$n >0$, the subgroup $_n A$ is finite. For $A$ a topological abelian group the notation $A^D$ will stand for the group of continuous homomorphisms $A\to\Q/\Z$.  The functor $A \mapsto A^D$
is an anti-equivalence of categories between $\ell$-primary torsion
abelian groups of cofinite type (with the discrete topology) and $\Z_{\ell}$-modules of finite type (with the profinite topology).

The notation $A^\wedge$ will stand for the profinite completion of $A$, i.e. the inverse limit of its finite quotients. We shall denote by $A_\wedge$ the inverse limit of the quotients $A/nA$ for all $n>0$, and by $A^{(\ell)}$ the inverse limit of the quotients
$A/\ell ^mA$ for a fixed prime $\ell$. Note that if $A$ is $\ell$-primary torsion of cofinite type, the map $A\to \varprojlim A/\ell^mA$ is surjective, inducing an isomorphism of finite groups $\overline A\stackrel\sim\to A^{(\ell)}$.

Unless otherwise stated, all cohomology groups will be taken with respect to the \'etale topology.\smallskip

\noindent {\em Tori.\/} For a torus $T$ over a field $F$ we shall denote by $\widehat T$ its character module and by $\widecheck T$ its module of cocharacters. These are finitely generated free $\Z$-modules equipped with a Galois action, and moreover $\widecheck T$ is the $\Z$-linear dual of $\widehat T$. The {\em dual torus} $T'$ of $T$ is by definition the $F$-torus with character group $\widecheck T$. The torus $T$ is {\em quasi-trivial} if $\widehat T$ is a permutation module, i.e. it has a Galois invariant basis. \smallskip

\noindent{\em Motivic complexes.\/} For $i\geq 0$ and a separated scheme $V$ of finite type over a field $F$ we denote by $z^i(V,\bullet)$ Bloch's cycle complex defined in \cite{hcg}. If $V$ is smooth, the \'etale motivic complex $\Z(i)$ over $V$ is defined as the complex of sheaves $z^i(-,\bullet)[-2i]$ on the small \'etale site of $V$. The analogously defined Zariski motivic complex is denoted by $\Z(i)_{\rm Zar}$. For an abelian group $A$ denote by $A(i)$ the complex $A\otimes\Z(i)$ and similarly for $\Z(i)_{\rm Zar}$ (since the terms of $\Z(i)$ are torsion free, this is the same as the derived tensor product). For $m$ invertible in $F$ we have a quasi-isomorphism of complexes of \'etale sheaves
\begin{equation}\label{Z(i)modm}
\Z/m\Z(i)\stackrel\sim\to\mu_m^{\otimes i}
\end{equation} where $\mu_m$ is the \'etale sheaf of $m$-th roots of unity placed in degree 0 (Geisser--Levine \cite{gl}, Theorem 1.5). Thus we shall also use the notation $\Q/\Z(i)$ for the direct limit of the sheaves $\mu_m^{\otimes i}$ for all $m>0$. \smallskip

\noindent{\em Function fields.\/} Throughout the paper, $k$ will be a finite extension of $\Q_p$ for a prime number $p$, and $X$ a smooth proper geometrically integral curve over $k$. The set of all closed points of $X$ will be denoted by $X^{(1)}$, and its function field by $K$. For a closed point $v\in X^{(1)}$ we denote by $\kappa(v)$ its residue field, and by $K_v$ the completion of $K$ for the discrete valuation induced by $v$. Then $\kappa(v)$ is also the residue field of the ring of integers $\calo_v$ and is a finite extension of $k$. Therefore $K_v$ is a {\em 2-dimensional local field\/}, i.e. a field complete with respect to a discrete valuation whose residue field is a classical local field.

\section{An Artin--Verdier style duality theorem} \label{one}

Let $T$ be a torus over $K$, and let $U\subset X$ be an open subscheme such that $T$ extends to a torus $\calt$ over $U$. The dual torus $T'$ then also extends to a torus $\calt'$ over $U$. Our first aim is to construct a duality pairing
\begin{equation}\label{caltpairing}
H^i(U, \calt)\otimes H^{3-i}_c(U, \calt')\to\Q/\Z.
\end{equation}

As above, denote by $\widehat\calt={\it Hom}(\calt, \G)$ the character group of $\calt$ and by $\widecheck\calt ={\it Hom}(\G, \calt)$ its cocharacter group. We have an quasi-isomorphism in the bounded derived category of \'etale sheaves \mbox{on $U$}
\begin{equation}\label{tordef}
\calt\cong \widecheck\calt\otimes^{\bf L} \Z(1)[1]
\end{equation}
using the fact that $\Z(1)[1]\cong\G$, and similarly,
\begin{equation}\label{dualtordef}
\calt'= \widehat\calt\otimes^{\bf L} \Z(1)[1].
\end{equation}
Composition of characters with cocharacters
$$
\widecheck \calt\otimes\widehat \calt\to{\it Hom}(\G,\G)\cong\Z
$$
together with the pairing of motivic complexes
$$\Z(1)\otimes^{\bf L}\Z(1)\to \Z(2)$$
(see e.g. \cite{totaro}) induces a pairing
\begin{equation}\label{pairing}
\calt\otimes^{\bf L} \calt'=(\widecheck \calt\otimes^{\bf L} \Z(1)[1])\otimes^{\bf L} (\widehat \calt\otimes^{\bf L} \Z(1)[1])\to \Z(2)[2].
\end{equation}

By passing to cohomology, we obtain pairings
$$
H^i(U, \calt)\otimes H^{3-i}_c(U, \calt')\to H^5_c(U, \Z(2)).
$$
Thus in order to obtain a pairing as in (\ref{caltpairing}), it remains to prove the following lemma which is implicit in the results of \cite{geisser} but not stated there.

\begin{lem}\label{globaltrace}  There is a canonical isomorphism $$H^5_c(U, \Z(2))\cong\Q/\Z.$$
\end{lem}

\begin{dem} Assume first $U=X$. Tensoring the exact sequence of abelian groups
$$
0\to \Z\to\Q\to\Q/\Z\to 0
$$
by $\Z(i)$ and using (\ref{Z(i)modm}) for all $m>0$ we obtain an exact sequence of complexes
\begin{equation}\label{Z(i)ex}
0\to \Z(i)\to\Q(i)\to \Q/\Z(i)\to 0
\end{equation}
Taking the associated long exact sequence for $i=2$ we obtain an exact sequence
$$
H^4(X, \Q(2))\to H^4(X, \Q/\Z(2))\to H^5(X, \Z(2))\to H^5(X, \Q(2)).
$$
Now denote by $\alpha$ the projection map from the small \'etale site of $X$ to its small Zariski site. According to (\cite{kahn}, Theorem 2.6. c)) the natural map $\Q(2)_{\rm Zar}\to \RR\alpha_*\Q(2)$ is an isomorphism in the derived category of Zariski sheaves. Hence there are isomorphisms
$$
H^j(X,\Q(2))\stackrel\sim\to H^j_{\rm Zar}(X,\Q(2))
$$
for all $j>0$. As $X$ has dimension 1 and $\Q(2)_{\rm Zar}$ is concentrated in degrees $\leq 2$ (see e.g. \cite{kahn}, Lemme 2.5), the latter are trivial for $j>3$, whence an isomorphism
$$
H^4(X, \Q/\Z(2))\stackrel\sim\to H^5(X, \Z(2)).
$$
By (\ref{Z(i)modm}) the group $H^4(X, \Q/\Z(2))$ is the direct limit of the groups $H^4(X, \mu_m^{\otimes 2})$ for all $m$. Setting $\Xbar:=X\times_k\kbar$, the degeneration of the Hochschild--Serre spectral sequence in \'etale cohomology together with the trace maps of Poincar\'e duality for $\Xbar$ and Tate duality for $k$ yield a chain of isomorphisms
$$
H^4(X, \mu_m^{\otimes 2})\cong H^2(k, H^2(\Xbar, \mu_m^{\otimes 2}))\cong H^2(k, \mu_m)\cong \Z/m\Z,
$$
whence the result by passing to the limit.

For $U\neq X$ consider the localization exact sequence
\begin{align*}
\bigoplus_{v\in X\setminus U} H^4(\kappa(v),i_v^*(\Z(2))) &\to H^5_c(U, \Z(2))\to \\ &\to H^5(X, \Z(2))\to \bigoplus_{v\in X\setminus U} H^5(\kappa(v), i_v^*(\Z(2)))
\end{align*}
constructed as in Proposition \ref{commute} (2) below (here $i_v$ denotes the closed immersion $\spec\kappa(v)\to X$). The groups $H^j(\kappa(v), i_v^*(\Z(2)))$ are 0 for $j>3$, since $\kappa(v)$ has cohomological dimension 2 and $i_v^*(\Z(2))$ is concentrated in degrees $\leq 2$. Hence the isomorphism ${H^5_c(U, \Z(2))\cong\Q/\Z}$ follows from the case $U=X$ in view of the localization exact sequence above.
\end{dem}

\begin{rema}\rm By similar arguments one can show that for $V$ smooth of dimension $d$ over a $p$-adic field $k$ there is an isomorphism $H^{2d+3}_c(V, \Z(d+1))\cong\Q/\Z.$

\end{rema}

We now prove the main theorem of this section.

\begin{theo} \label{artverd}
 The pairing (\ref{caltpairing})
induces a perfect pairing of finite groups
$$(H^i(U,\T)\{\ell\})^{(\ell)}  \times H^{3-i}_c(U,\T')^{(\ell)}\{ \ell \}
 \to \Q_{\ell}/\Z_{\ell}$$
for each prime number $\ell$.
\end{theo}

\dem This is essentially the same
argument as in (\cite{dhsza1}, Theorem 3.4) or (\cite{adt}, Theorem II.5.2). Recall that given a finite group scheme $\calf$ over $U$, combining Poincar\'e duality for $U\times_k\kbar$ with Tate--Nakayama local duality over $k$ yields a perfect pairing
of finite groups
\begin{equation} \label{finidual}
H^i(U,\calf) \times H^{4-i}_c(U,\calf') \to \Q/\Z
\end{equation}
where $\calf':=\Hom(\calf,\Q/\Z(2))$ (see e.g. \cite{geisser}, \S 5.2). We shall apply this with $\calf={}_m\T$. Note that in this case $\calf'$ is indeed the $m$-torsion in $\calt'$, since in view of (\ref{tordef}) we have
$$
\Hom({}_m\T,\mu_m^{\otimes 2})\cong \Hom(\widecheck \T\otimes\mu_m,\mu_m^{\otimes 2})\cong \Hom(\widecheck \T,\mu_m)\cong {}_m\T'.
$$
Using the Kummer sequence
$$0 \to H^{i-1}(U,\T)/\ell ^n \to H^i(U, {}_{\ell^n}\T) \to
{}_{\ell^n} H^i(U,\T) \to 0$$
and passing to the direct limit, we get an exact sequence
$$0 \to H^{i-1}(U,\T)\otimes {\Q_{\ell}/\Z_{\ell}}
\to \varinjlim_n H^i(U, {}_{\ell^n}\T) \to
H^i(U,\T) \{ \ell \} \to 0$$
Since $H^{i-1}(U,\T) \otimes \Q_{\ell}/\Z_{\ell}$ is divisible, this induces
for each $m >0$ an isomorphism
$$(\varinjlim_n H^i(U, {}_{\ell^n} \T))/\ell^m
\simeq (H^i(U,\T) \{ \ell \})/\ell^m.$$
Taking the inverse limit over $m$, we obtain
an isomorphism
$$(\varinjlim_n H^i(U, {}_{\ell^n}\T))^{(\ell)} \simeq
(H^i(U,\T) \{ \ell \})^{(\ell)}.$$

\noindent Similarly, the exact sequence
$$0 \to H^{3-i}_c(U,\T')/\ell^n
\to H^{4-i}_c(U, {}_{\ell^n}\T') \to
{}_{\ell^n}H^{4-i}_c(U,\T')  \to 0$$
 yields after passing to the inverse limit an isomorphism
$$(\varprojlim_n H^{4-i}_c(U,{}_{\ell^n}\T'))\{ \ell \} \simeq
H^{3-i}_c(U,\T')^{(\ell)} \{\ell\}$$
because the $\ell$-adic Tate module of $H^{4-i}_c(U,\T')$ is torsion free and $H^{4-i}_c(U, {}_{\ell^n}\T')$ is finite.
Now the result follows from the duality (\ref{finidual})
applied to each ${}_{\ell^n}\T$.
\enddem

\begin{rema}\label{avrema}\rm The theorem could also be derived from the stronger result of (\cite{schvh}, Theorem 4.8) which  works for arbitrary $U$ and not only for $U$ sufficiently small. In both approaches the argument is based on the perfectness of the duality pairing (\ref{finidual}) for finite modules. Note, however, that in \cite{schvh} the pairing for tori itself is constructed by patching together duality pairings for finite modules, whereas we have worked with the canonically defined pairing (\ref{caltpairing}).
\end{rema}

\section{Local duality}

We now handle the local theory over the 2-dimensional local field $K_v$. By a similar argument as in Lemma \ref{globaltrace}, we prove:

\begin{lem} For each closed point $v\in X$ there is a canonical isomorphism $$H^4(K_v, \Z(2))\cong\Q/\Z.$$
\end{lem}

\begin{dem} Exact sequence (\ref{Z(i)ex}) over the spectrum of $K_v$ yields a long exact sequence
$$
H^3(K_v, \Q(2))\to H^3(K_v, \Q/\Z(2))\to H^4(K_v, \Z(2))\to H^4(K_v, \Q(2))
$$
As $\spec(K_v)$ has dimension 0, the arguments of the proof of Lemma \ref{globaltrace} yield the vanishing of $H^j(K_v, \Q(2))$ for $j>2$. Thus the map
$$
H^3(K_v, \Q/\Z(2))\to H^4(K_v, \Z(2))
$$
is an isomorphism. Taking isomorphism (\ref{Z(i)modm}) again into account, we have
\begin{equation}\label{kato}
H^3(K_v, \Q/\Z(2))\cong\Q/\Z
\end{equation} by Kato's class field theory for higher local fields \cite{kato} (actually, this case is easy; see also \cite{adt}, Theorem 2.17).
\end{dem}

Given a torus $T$ over $K_v$, we obtain duality pairings
\begin{equation}\label{localpairing}
H^i(K_v, T)\otimes H^{2-i}(K_v, T')\to H^4(K_v, \Z(2))\cong\Q/\Z
\end{equation}
by similar arguments as during the construction of the pairing (\ref{caltpairing}) in the previous section.

Recall also (e.g. from \cite{adt}, I. 2.17) that for a finite $\gal({\overline K}_v|K_v)$-module $F$
there are perfect pairings of finite groups
\begin{equation}\label{finitepairing}
H^i(K_v, F)\otimes H^{3-i}(K_v, F')\to H^3(K_v, \Q/\Z(2))\cong\Q/\Z.
\end{equation}
where $F':=\Hom(F,\Q/\Z(2))$.

Using this result, we prove :

\begin{prop} \label{localdual}
The pairing (\ref{localpairing}) is a perfect pairing of finite groups for $i=1$. For $i=0$ it becomes a perfect pairing between a profinite and a torsion group, after replacing $H^0(K_v, T)$ by its profinite completion.
\end{prop}

\begin{dem}
Given $n>0$, the exact sequences
$$
0\to {}_nT\to T\to T\to 0
$$
and
$$
0\to {}_nT'\to T'\to T'\to 0
$$
induce a commutative exact diagram
$$
\begin{CD}
0 @>>> H^0(K_v, T)/n@>>> H^1(K_v, {}_nT)@>>> {}_nH^1(K_v, T)@>>> 0\\
&& @VVV @VV{\cong}V @VVV \\
0 @>>> ({}_nH^2(K_v, T'))^D@>>> H^2(K_v, {}_nT')^D @>>> (H^1(K_v, T')/n)^D @>>> 0
\end{CD}
$$
Here the middle isomorphism is induced by the perfect pairing (\ref{finitepairing}), taking into account that ${}_nT'\cong \Hom({}_nT,\Q/\Z(2))$ as in the previous section.  The groups $H^1(K_v,T)$ and $H^1(K_v, T')$ are finite because they have finite exponent as a consequence of Hilbert's Theorem 90 and they are of finite cotype because the groups $H^1(K_v, {}_nT)$ and $H^1(K_v, {}_nT')$ are finite as recalled above (see e.g. \cite{cogal}, III.4.3 for details). Therefore by taking $n$ sufficiently large we see that the map $H^1(K_v,T)\to H^1(K_v, T')^D$ is surjective; in particular, the order of $H^1(K_v,T)$ is greater than that of $H^1(K_v, T')$. Exchanging the roles of $T$ and $T'$ we obtain the statement for $i=1$.

Once we know the case $i=1$, we also know that the right vertical map in the diagram above is an isomorphism. But then so is the left vertical map and we obtain the case $i=0$ by passing to the inverse limit over all $n$.
\end{dem}

\begin{rema}\rm We also have $H^i(K_v, T)=0$ for $i\geq 3$. For $i>3$ this follows from the fact that $K_v$ has cohomological dimension 3, and for $i=3$ one can prove this by a similar argument as in (\cite{cogal}, II.5.3).
\end{rema}

\section{Global duality: preliminary results}

This section is devoted to basic facts needed for the proof of the global duality theorem. We begin with some general statements on compact support cohomology which are in fact valid over an arbitrary base field.

\begin{prop}\label{commute} Let $U\subset X$ be an open subscheme and $\calg$ an \'etale sheaf on $U$.

\noindent $(1)$ We have a long exact sequence of cohomology groups
$$
\cdots\to H^i_c(U,\calg)\to H^i(U, \calg)\to\bigoplus_{v\in X\setminus U} H^i(K_v^h, \calg_v)\to H^{i+1}_c(U, \calg)\to\cdots
$$
where $K_v^h$ is the henselisation of $K$ with respect to the discrete valuation induced by $v$ and $\calg_v$ is the pullback of $\calg$ by the natural map $\spec K_v^h\to U$.

\noindent $(2)$ If $V\subset U$ is a further open subscheme, we have a long exact sequence
$$
\cdots\to H^i_c(V,\calg)\to H^i_c(U, \calg)\to\bigoplus_{v\in U\setminus V} H^i(\kappa(v), i_v^*\calg)\to H^{i+1}_c(V, \calg)\to\cdots
$$
where $i_v$ is the inclusion map of the closed point $v$ in $U$.

\noindent $(3)$ {\rm (Three Arrows Lemma)} The diagram
$$
\begin{CD}
H^i_c(V,\calg) @>>> H^i_c(U,\calg) \cr
@VVV @VVV \cr
H^i(V,\calg) @<<< H^i(U,\calg)
\end{CD}
$$
commutes.
\end{prop}

\begin{dem} $(1)$ is proven by exactly the same argument as (\cite{adt}, Lemma II.2.4).
To prove (2), we use the exact sequence
$$
\cdots\to H^i_c(V,\calg)\to H^i_c(U, \calg)\to H^i_c(U\setminus V, i^*\calg)\to\cdots
$$
of (\cite{milne}, III.1.30) for the closed immersion $i:\, (U\setminus V)\to U$, and note that
$$
H^i_c(U\setminus V, i^*\calg)= H^i(U\setminus V, i^*\calg)\cong \bigoplus_{v\in U\setminus V}H^i(\kappa(v), i_v^*\calg)
$$
as $U\setminus V$ is a finite set of closed points.

For (3), denote by $j$ the inclusion map $U\to X$ and by $j_U$ the inclusion map $V\to U$. Since $H^i_c(U,\calg)=H^i(X, j_!\calg)$ and similarly for $V$, we can rewrite the diagram in terms of Ext-groups of sheaves as
$$
\begin{CD}
\Ext^i_X(\Z,(j\circ j_U)_!j^*_U\calg) @>>> \Ext^i_X(\Z,j_!\calg) \cr
@VVV @VVV \cr
\Ext^i_V((j\circ j_U)^*\Z,(j\circ j_U)^*(j\circ j_U)_!j^*_U\calg) @<<< \Ext^i_U(j^*\Z,j^*j_!\calg)
\end{CD}
$$
where the vertical maps are induced by composing with $(j\circ j_U)^*$ and $j^*$, respectively, and the upper horizontal map is induced by the natural map $j_{U!}j^*_U\calg\to\calg$. To recognize the lower horizontal map as a restriction map, use the identities $j^*j_!=\id_U$ and $j_U^*j_{U!}=\id_V$ to rewrite the diagram as
$$
\begin{CD}
\Ext^i_X(\Z,(j\circ j_U)_!j^*_U\calg) @>>> \Ext^i_X(\Z,j_!\calg) \cr
@VVV @VVV \cr
\Ext^i_V(j_U^*j^*\Z,j^*_U\calg) @<<< \Ext^i_U(j^*\Z,\calg)
\end{CD}
$$
where the lower horizontal map is induced by $j_U^*$.
To verify the commutativity of this diagram, note that the left vertical map factors as
$$
\Ext^i_X(\Z,(j\circ j_U)_!j^*_U\calg)\to \Ext^i_U(j^*\Z,j_{U!}j^*_U\calg)\to \Ext^i_V(j^*_Uj^*\Z, j^*_U\calg)
$$
where the second map factors through $\Ext^i_U(j^*\Z,\calg)$ in view of the commutative diagram
$$
\begin{CD}
\Ext^i_U(j^*\Z, j_{U!}j^*_U\calg) @>>> \Ext^i_U(j^*\Z,\calg) \cr
@VVV @VVV \cr
\Ext^i_V(j_U^*j^*\Z,j^*_Uj_{U!}j^*_U\calg) @>>> \Ext^i_V(j_U^*j^*\Z,j_U^*\calg)
\end{CD}
$$
in which the horizontal maps are induced by the map $j_{U!}j^*_U\calg\to\calg$ and the vertical maps by $j_U^*$. The desired commutativity thus reduces to commutativity of the diagram
$$
\begin{CD}
\Ext^i_X(\Z,(j\circ j_U)_!j^*_U\calg) @>>> \Ext^i_X(\Z,j_!\calg) \cr
@VVV @VVV \cr
\Ext^i_U(j^*\Z,j_{U!}j^*_U\calg) @>>> \Ext^i_U(j^*\Z,\calg)
\end{CD}
$$
which holds by functoriality of the adjunction map for $j_!$ with respect to the map $j_{U!}j^*_U\calg\to\calg$.
\end{dem}

Let now $T$ be a $K$-torus with dual torus $T'$. We shall work with
nonempty Zariski open subsets $U$ of $X$ which will always be assumed
affine and such that $T$, $T'$ extend to $U$-tori $\T$, $\T '$.

The next two corollaries are in fact valid over any base field of characteristic 0.

\begin{cor}\label{hensel} The natural maps
$$
H^i(K_v^h, T)\to H^i(K_v,T)
$$
are isomorphisms for $i\geq 1$. Consequently, we have a long exact sequence
\begin{equation}\label{locseq}
\cdots\to H^i_c(U,\calt)\to H^i(U, \calt)\to\bigoplus_{v\in X\setminus U} H^i(K_v, T)\to H^{i+1}_c(U, \calt)\to\cdots
\end{equation}
starting from degree 1.
\end{cor}

\begin{dem} The first statement is a consequence of Greenberg's approximation theorem \cite{greenberg} (see e.g. the proof of \cite{dhsza1}, Lemma 2.7 for a detailed argument). The second one is a consequence of Proposition \ref{commute} (1).
\end{dem}

\begin{cor}
The groups $H^i(U,\T)$ and $H^i_c(U,\T)$ are torsion for $i\geq 2$.
\end{cor}

\begin{dem} Using the exact sequence of the previous corollary and the fact that the groups $H^i(K_v, T)$ are torsion for $i>0$ we see that it is enough to treat the groups $H^i(U,\T)$. Next, trivializing $\T$ by means of a finite \'etale cover of $U$ and using a restriction-corestriction argument we reduce to the case $\T=\G^r$, and finally to $T=\G$ by compatibility of cohomology with finite products. The case $i=2$ then follows from the fact that $H^2(U,\G)=\br U$ is a subgroup of the torsion group $\br K$. For $i>2$ we use the isomorphism $H^i(U,\G)\cong H^{i+1}(U,\Z(1))$ and the piece of the long exact sequence
$$
H^{i}(U,\Q/\Z(1))\to H^{i+1}(U,\Z(1))\to H^{i+1}(U,\Q(1)).
$$
As $H^{i+1}(U,\Q(1))=0$ for $i\geq 2$ by a similar argument as in the proof of Lemma \ref{globaltrace}, the first map is surjective for $i\geq 2$, but $H^{i}(U,\Q/\Z(1))$ is torsion.
\end{dem}

We now return to our basic assumption that the base field $k$ is a finite extension of $\Q_p$. For each $i \geq 1$, we set
$$\Sha^i(T):=\ker (H^i(K,T) \to \prod_{v \in X^{(1)}} H^i(K_v,T)).$$

\begin{prop} \label{firstprop} ${}$

\noindent $(1)$ The group
$H^1(U,\T)_{\rm tors}$ is of cofinite type, and so are the torsion groups $H^i(U,\T)$
and $H^i_c(U,\T)$ for $i \geq 2$.

\noindent $(2)$ The group $H^1(K,T)$ is of finite exponent
and the groups $\Sha^i(T)$ are finite for $i=1,2$.
\end{prop}

\begin{dem} (1) For each $n>0$  the Kummer sequence yields a surjection $H^i(U,{}_n\T) \twoheadrightarrow
{}_n H^i(U,\T)$, and we have already recalled during the proof of Theorem \ref{artverd} that the groups  $H^i(U,{}_n\T)$ are finite. Thus $H^i(U,\T)_{\rm tors}$ is of cofinite type (and equals $H^i(U,\T)$ for $i\geq 2$).  To see that $H^i_c(U,\T)$ is of cofinite type for $i\geq 2$,
we use the exact sequence
$$\bigoplus_{v \not \in U} H^{i-1}(K_v,T) \to H^i_c(U,\T) \to H^i(U,T)$$
coming from (\ref{locseq}) and the fact that the groups $H^{i-1}(K_v,T)$ are torsion
and of cofinite type (the latter because for every $n >0$, the
finite group $H^{i-1}(K_v,{}_n T)$ surjects
onto ${}_n H^{i-1}(K_v,T)$).

\smallskip

(2) If $L|K$ is a finite extension of degree $d$ that splits $T$, we have  $H^1(L,\G)=0$ by Hilbert's Theorem 90, whence $d\, H^1(K,T)=0$ by a restriction-corestriction argument. To see the finiteness of $\Sha^1(T)$, note first that by exact sequence (\ref{locseq}) the intersection in $H^1(K,T)$ of the images of the composite maps $H^1_c(U,\T)\to H^1(U,\T)\to H^1(K,T)$ for all $U$ is exactly $\Sha^1(T)$, so in particular $\Sha^1(T)$ is contained in the image of the map $H^1(U, \T)\to H^1(K,T)$ for a fixed $U$. But this map factors through $H^1(U, \T)/d$ since $d\, H^1(K,T)=0$. On the other hand, the group  $H^1(U, \T)/d$ injects in $H^2(U, {}_d\T)$ by the Kummer sequence. As we have just recalled, the latter group is finite.

In order
to prove that $\Sha^2(T)$ is of finite exponent it suffices
to prove ${\Sha^2(\G)=0}$ by another restriction-corestriction argument. Pick $\alpha \in \Sha^2(\G) \subset
H^2(K,\G)=\br K$. Since $\alpha$ maps to 0 by the localisation maps $\br K\to \br K_v$, it also maps to 0 by the residue maps $\br K\to H^1(\kappa(v),\Q/\Z)$, which means that $\alpha \in \br X$. In particular $\alpha\in\br\calo_{X,v}$, where $\calo_{X,v}$ is the local ring of $X$ at the closed point $v$ whose completion is the ring of integers $\calo_v$ of $K_v$. Since the  restriction map $\br\calo_v\to\br K_v$ is injective, we conclude that $\alpha$ maps to 0 already in $\br\calo_v$. Composing with the isomorphism $\br\calo_v\cong\br\kappa(v)$ we conclude that $\alpha$ maps to 0 by the `evaluation at $v$' map $\br X\to\br\kappa(v)$ induced by $v$.
Therefore $\alpha$ is orthogonal to $\pic X$ for Lichtenbaum's duality
pairing $\br X\times \pic X\to\Q/\Z$ which is defined in \cite{licht} by a sum of evaluation maps, from which we conclude $\alpha=0$.

\smallskip

It remains to prove that $\Sha^2(T)$ is of cofinite type. We observe that
each element $\alpha$ of $\Sha^2(T) \subset H^2(K,T)$ is in the image of
$H^2(V,\T)$ for some open subset $V \subset U$. Exact sequence
(\ref{locseq}) shows that $\alpha$ comes from $H^2_c(V,\T)$, hence also from $H^2_c(U,\T)$ by the Three Arrows Lemma (Proposition \ref{commute} (3)). This shows that
$\Sha^2(T)$ is a subquotient of $H^2_c(U,\T)$ (for a fixed $U$), hence
it is of cofinite type by statement (1).
\end{dem}

\begin{rema}\label{shaq}\rm We note for further reference that the proof of part (2) above also yields the triviality of $\Sha^2(Q)$ for a quasi-trivial torus $Q$. Indeed, by Shapiro's lemma it suffices to establish the case $Q=\G$ which was treated above.

\end{rema}

Next, define $$\D^2(U,\T)= \im [H^2_c(U,\T) \to H^2(K,T)].$$ We observe
that if $V \subset U$, then $\D^2(V,\T) \subset \D^2(U,\T)$ by
covariant functoriality of the groups $H^2_c$.

 The following statement will be a key tool in the proof of the global duality theorem.

\begin{prop} \label{station}
For a fixed prime number $\ell$ there exists a nonempty Zariski
open subset $U_0\subset X$ such that for every nonempty Zariski open
 $U \subset U_0$ we have
$$\D^2(U,\T)\{ \ell \}=\D^2(U_0,\T) \{ \ell \}=\Sha^2(T)\{ \ell \}; $$
consequently, the first two groups are finite.
\end{prop}

The proof uses the following lemma which is probably well known.

\begin{lem} \label{stationl}
Every decreasing sequence $(A_n)$ of abelian $\ell$-primary torsion groups
of cofinite type stabilizes.
\end{lem}

\dem We shall prove the dual statement: every inductive system $(B_n)$ of finitely generated $\Z_{\ell}$-modules with surjective transition maps stabilizes. Since $\Z_{\ell}$
is a principal ideal domain, a finitely generated $\Z_{\ell}$-module $B$ is a direct sum of the finite submodule $B_{\rm tors}$  and of a free module of finite rank equal to the $\Q_\ell$-dimension of $B\otimes_{\Z_\ell}\Q_\ell$. The surjectivity of the transition maps implies that the rank
of $B_n$ is non-increasing,  hence
we can assume that the sequence $(B_n)$ has constant rank. Also,
the kernel of every map $B_n \to B_{n+1}$ is finite, which implies
that the induced map $(B_n)_{\rm tors} \to (B_{n+1})_{\rm tors}$ is
surjective. Since each group $(B_n)_{\rm tors}$ is finite, this implies
that the sequence $(B_n)_{\rm tors}$ stabilizes, whence the lemma.
\enddem

\noindent{\em Proof of Proposition \ref{station}.} Fix an open subset $V_0$ of $X$ such that $T$
extends to a torus over $V_0$. For open subsets $V \subset U$ of
$V_0$, we observe that $\D^2(U,\T)\{ \ell \}$ is of cofinite type
by Proposition~\ref{firstprop},
and $\D^2(V,\T)\{\ell\} \subset \D^2(U,\T)\{\ell\}$. By
Lemma~\ref{stationl}, there exists a nonempty Zariski
open subset $U_0$ of $X$ such that for every nonempty Zariski open
subset $U \subset U_0$, we have $\D^2(U,\T)\{ \ell \}=\D^2(U_0,\T) \{ \ell \}$.
Since we have
$$\D^2(U,\T) \subset \ker [H^2(K,T)
\to \prod_{v \not \in U} H^2(K_v,T)] ,$$
we obtain $$\D^2(U_0,\T) \{ \ell \}=\Sha^2(T)\{ \ell \} .$$
By Proposition~\ref{firstprop}, (2), we get that
for $U \subset U_0$, the group $\D^2(U,\T)\{ \ell \}$ is finite.
\enddem

\section{The global duality theorem}

We keep notation from the previous sections. This section is devoted to the proof of the following theorem.

\begin{theo}\label{dualsha} Given a $K$-torus $T$, there is a perfect pairing of finite
groups
$$\Sha^1(T) \times \Sha^2(T') \to \Q/\Z.$$
\end{theo}

Besides the results of the previous sections, we shall need one auxiliary statement. To state it, we first construct a map
\begin{equation}\label{map}
\bigoplus_{v \in X^{(1)}} H^1(K_v,T) \to H^2_c(U,\T).
\end{equation}

Each element $\alpha\in\bigoplus_{v \in X^{(1)}} H^1(K_v,T)$ comes from $\bigoplus_{v \not \in V}H^1(K_v,T)$ for some  nonempty open subset
$V \subset U$. Using the exact sequence of Corollary \ref{hensel}  for $i=2$ and $\calg=\T$
 we can thus send $\alpha$ to $H^2_c(V,\T)$ and, by covariant functoriality of compact support cohomology, to  $H^2_c(U,\T).$ This map does not depend on the choice of $V$ as for $W\subset V$ the diagram
$$
 \begin{CD}
\bigoplus_{v \not \in W} H^1(K_v,T) @>>> H^2_c(W,\T)  \cr
 @AAA @VVV  \cr
\bigoplus_{v \not \in V} H^1(K_v,T) @>>>  H^2_c(V,\T)
\end{CD}
$$
commutes.

Now we can state:

\begin{prop} \label{exactd}
The sequence
\begin{equation} \label{exactall}
\bigoplus_{v \in X^{(1)}} H^1(K_v,T) \to H^2_c(U,\T) \to \D^2(U,\T) \to 0
\end{equation}
is exact.
\end{prop}

\noindent{\em Proof.} Let $V \subset U$ be a Zariski open subset. By the Three Arrows Lemma (Proposition~\ref{commute} (3)) and  Corollary \ref{hensel}
there is a
commutative diagram with exact upper row
\begin{equation}\label{diag}
\begin{CD}
\bigoplus_{v \not \in V} H^1(K_v,T) @>>> H^2_c(V,\T) @>>> H^2(V,\T) \cr
&& @VVV @AAA  \cr
&&  H^2_c(U,\T) @>>> H^2(U,\T).
\end{CD}
\end{equation}

\noindent
Given $\alpha \in \oplus_{v \in X^{(1)}} H^1(K_v,T)$, choose $V$ as in the construction of the map (\ref{map}). By exactness of the upper row the image of $\alpha$ in
$H^2(V,\T)$ is zero, hence so is its image in $H^2(K,T)$. In other words, the image of $\alpha$
in $H^2_c(V,\T)$ maps to 0 in $H^2(K,T)$. By commutativity of the square in the diagram so does its image in $H^2_c(U,T)$, which proves that the sequence (\ref{exactall}) is a complex.

Conversely, pick $\beta \in \ker [H^2_c(U,\T) \to \D^2(U,\T)]$ and a nonempty open subset $V\subset U$.
 Consider the commutative diagram
\begin{equation}\label{diag2}
\begin{CD}
H^2_c(V,\T) @>>> H^2_c(U,\T) @>>> \bigoplus_{v\in U\setminus V} H^2(\kappa(v),i_v^* \T) \\
&& @VVV @VVV \\
&& H^2(K,T) @>>> \bigoplus_{v\in U\setminus V} H^2(K_v,T)
\end{CD}
\end{equation}
in which the upper row comes from the exact sequence of Proposition \ref{commute} (2) and the left vertical map is given by the composite map $ H^2_c(U,\T)\to  H^2(U,\T)\to  H^2(K,T)$. The right vertical map can be constructed as a composite
$$
H^2(\kappa(v),i_v^* \T)\stackrel\sim\leftarrow H^2(\calo_{U,v}^h, \T)\to H^2(K_{v}^h, T)\stackrel\sim\to H^2(K_v,T)
$$
where the first isomorphism holds by (\cite{adt}, Proposition II.1.1 (b)).
In fact, this map is injective. To see this, it suffices to show injectivity of the middle arrow, which may be identified with a map $H^2(G/I,T(\overline K)^I)\to H^2(G,T(\overline K))$ coming from a Hochschild--Serre spectral
sequence, where $I$ is the inertia subgroup in $G=\gal(\overline K|K_v^h)$.  The latter map is injective because $H^0(G/I, H^1(I,T(\overline K)))=0$ by Hilbert's Theorem 90 as $T$ is split by an unramified extension.

Commutativity of the square in the diagram follows from the commutativity of
$$
\begin{CD}
H^2_c(U,\T) @>>> \bigoplus_{v\in U\setminus V} H^2(\kappa(v),i_v^* \T) \\
@VVV @AA{\cong}A \\
H^2(U,\T) @>>> \bigoplus_{v\in U\setminus V} H^2(\calo_{U,v}^h,\T)
\end{CD}
$$
which is just a diagram of restriction maps.

As the image of our element $\beta$ in $H^2(K, T)$ is zero, diagram (\ref{diag2}) shows that its image in each $H^2(\kappa(v),i_v^* \T)$ is also zero, and therefore $\beta$ comes from an element $\gamma\in H^2_c(V,\T)$. Since the image of $\gamma$ in $H^2(K,T)$ is 0, we can choose $V$ so small that $\gamma$ maps to 0 already by the map $H^2_c(V,\T)\to H^2(V,\T)$.  The exact upper row of diagram (\ref{diag}) then yields the exactness of the sequence of the proposition.
\enddem

\noindent{\em Proof of Theorem \ref{dualsha}.} We already know from
Proposition~\ref{firstprop} (2) that $\Sha^1(T)$ and $\Sha^2(T')$ are finite. We construct a perfect pairing
$$\Sha^1(T)\{\ell\} \times \Sha^2(T')\{\ell\} \to \Q/\Z$$
for each prime number $\ell$. Define the group $D^1_{sh}(U,\T) $ by the exact sequence
$$0 \to D^1_{sh}(U,\T) \to H^1(U,\T) \to \prod_{v \in X^{(1)}} H^1(K_v,T).$$

As the groups $H^1(K_v,T)$ have exponent bounded by the degree of a field extension splitting $T$, their product has trivial maximal divisible subgroup, and the same holds for its $\ell$-primary torsion part.
Thus  the maximal divisible subgroup of
$H^1(U,\T)\{\ell\}$ is a subgroup of $D^1_{sh}(U,\T)\{\ell\}$. This fact, together with
Proposition~\ref{exactd} yields the exact rows of the commutative diagram
$$
\begin{CD}
0 @>>> \overline{D^1_{sh}(U,\T) \{\ell\}} @>>> \overline{H^1(U,\T)\{\ell\}}
@>>> \prod_{v \in X^{(1)}} H^1(K_v,T)\{ \ell \} \cr
&&&& @VVV @VV{\cong}V \cr
0 @>>> \D^2(U,\T')\{\ell\} ^D @>>> H^2_c(U,\T')\{ \ell \}^D @>>>
(\bigoplus_{v \in X^{(1)}} H^1(K_v,T')\{\ell\} )^D.
\end{CD}
$$
The right vertical map is the isomorphism given by local duality (Proposition~\ref{localdual}). The middle vertical map is induced by the duality pairing (\ref{caltpairing}), noting that the maximal divisible subgroup of $H^1(U,\T)\{\ell\}$ is annihilated by a pairing with a torsion group.  On the other hand, by
Proposition~\ref{station}
we have a nonempty Zariski open subset $U_0$ such that
for every $U \subset U_0$ we have the equalities
$$\D^2(U,\T')\{\ell \}=\D^2(U_0,\T')\{\ell\}=\Sha^2(T')\{\ell\} .$$

Moreover, Proposition \ref{exactd} implies that in this case the group $H^2_c(U,\T')\{\ell\}$ is finite since it is of cofinite type by Proposition~\ref{firstprop} (1) and the groups $H^1(K_v,T)$ have bounded exponent. Also, since $H^1(U,\T)\{\ell\}$ is of cofinite type by Proposition~\ref{firstprop} (1), we have $H^1(U,\T)\{\ell\}^{(\ell)}\cong\overline{ H^1(U,\T)\{\ell\}}$ (see our preliminary remarks on abelian groups). By Theorem~\ref{artverd}, the middle vertical map is therefore an isomorphism for
$U \subset U_0$.
Finally, the diagram yields for $U \subset U_0$ an isomorphism
$$\overline{D^1_{sh}(U,\T) \{\ell\}} \stackrel\sim\to \D^2(U,\T')\{\ell\} ^D =\Sha^2(T')\{\ell\} ^D.$$
Since ${{\Sha^1(T)} \{\ell\}}$ is the direct limit of the
${D^1_{sh}(U,\T) \{\ell\}}$ for $U \subset U_0$, we conclude
$${\Sha^1(T) \{\ell\}}\cong
\Sha^2(T')\{\ell\} ^D $$
since a
direct limit of divisible groups is divisible, and the proof is complete.
\enddem

\begin{rema}\label{exsha}\rm The theorem above is nonvacuous only if one knows that there exist tori over $K$ with nontrivial Tate--Shafarevich group. In other words, one needs examples of torsors under $K$-tori that have points everywhere locally but not globally. Such examples were constructed recently by Colliot-Th\'el\` ene, Parimala and Suresh (\cite{ctps2}, Proposition 5.9 and Remark 5.10).

Here is a simpler example, communicated to us by J-L. Colliot-Th\'el\`ene. Over $K=\Q_p(t)$ consider the torus $T$ of equation
$$
(x_1^2-ty_1^2)(x_2^2-(1+t)y_2^2)(x_3^2-t(1+t)y_3^2)=1
$$
and the $K$-torsor $Y$ of equation
$$
(x_1^2-ty_1^2)(x_2^2-(1+t)y_2^2)(x_3^2-t(1+t)y_3^2)=p
$$
under $T$. The elements $t$, $t+1$, $t(t+1)$ are units in the local ring $A$ of $\P^1_{\Z_p}$ at the generic point of the special fibre, and map to nonsquares in its residue field. It follows that every nonzero value of
$$
(x_1^2-ty_1^2)(x_2^2-(1+t)y_2^2)(x_3^2-t(1+t)y_3^2)
$$
with $x_i, y_i\in \Q_p(t)$ has even valuation in the discrete valuation $w$ of $A$, hence cannot equal $p$ which is a uniformizer of this valuation. Thus $Y$ has no points over $K$. On the other hand, assuming $p$ odd and $-1$ a square in $\Q_p$, one checks that for every closed point $v$ of $\P^1_{\Q_p}$ one of $t$, $t+1$ or $t(t+1)$ is a square in $K_v$, and therefore $Y(K_v)\neq\emptyset$.

The examples of \cite{ctps2} have the additional property that they have points over all completions of $K$ with respect to its discrete valuations (not just those coming from points of the curve). Here there is no point over $K_w$.
\end{rema}

We would also like to record here the following analogue of Theorem \ref{dualsha} for finite modules which is used in our paper \cite{pcwa}.

\begin{theo}\label{dualfinisha}
Let $F$ be a finite Galois module over $K$.
There are perfect pairings of finite
groups
$$\Sha^i(F) \times \Sha^{4-i}(F') \to \Q/\Z$$
for $i=1,2$.
\end{theo}

\begin{dem} The proof follows exactly the same pattern as in the case
of tori, so we only point out the changes to be made. Let $U$ be a nonempty Zariski open subset of $X$ such that $F$ extends to
an \'etale group scheme $\calf$ over $U$.
We have exact sequences

$$
\bigoplus_{v \in X^{(1)}} H^1(K_v,F) \to H^2_c(U,\calf) \to \D^2(U,\calf)
\to 0$$
and
$$\bigoplus_{v \in X^{(1)}} H^2(K_v,F) \to H^3_c(U,\calf) \to \D^3(U,\calf)
\to 0.$$

Indeed, the only difference with the proof of Proposition \ref{exactd} is that for
$v \in U$, the injectivity of the map $H^i(\calo_v,\calf) \simeq
H^i(\kappa(v),\calf) \to
H^i(K_v,F)$ for $i \geq 0$  is proven by the following argument. With notation as there, we have to prove the injectivity of the inflation map
$H^i(G_v/I_v,F(K_v ^{\rm nr}) \to H^i(G,F(\ov K_v))$. Since by assumption $F$ is unramified at $v$, we have an equality of Galois modules $F(\ov K_v)=F(K_v ^{\rm nr})$.
Therefore injectivity
follows from the fact that the canonical
surjection $G_v \to G_v/I_v$ has a section (\cite{cogal}, \S II.4.3).

The analogue of Corollary~\ref{station} is here immediate because we
already know that the groups $H^i_c(U,\calf)$ are finite for $i \geq 0$, hence
we get
$$\Sha^i(F)={\D}^i(U,\calf) , i=2,3$$
for $U$ sufficiently small. We also already know the analogue of
Theorem~\ref{artverd} without having to mod out by divisible parts, which
makes the end of the argument similar (but simpler), using
the duality between $H^i(U,\calf)$ and $H^{4-i}(U,\calf')$ for $i=1,2$ and
the local duality between $H^i(K_v,F)$ and $H^{3-i}(K_v,F')$ for
$i=1,2$.
\end{dem}

\section{A cohomological obstruction to the Hasse principle for tori} \label{five}

As before, let $K$ be the function field of a smooth projective geometrically
integral curve $X$ defined over a finite extension $k$ of $\Q_p$. Given a smooth geometrically integral $K$-scheme $Y$, for $U\subset X$ sufficiently small we find a smooth geometrically integral $U$-scheme $\mathcal Y$ such that ${\mathcal Y}\times_UK\cong Y$. Thus we may define
$$
Y(\A_K):=\varinjlim_{V\subset U} \prod_{v\not\in V}Y(K_v)\times  \prod_{v\in V}{\mathcal Y}(\calo_v)
$$
noting that $\calo_v$ is a completion of a local ring of the open subscheme $V\subset U$.

Consider now the group $H^3(Y, \Q/\Z(2))$. It is the direct limit of the groups $H^3({\mathcal Y}\times_UV,\Q/\Z(2))$ for $V\subset U$ open.
Recall also that $H^3(K_v,\Q/\Z(2))\cong\Q/\Z$ by (\ref{kato}) and
\begin{equation}\label{vanish}
H^3(\calo_v,\Q/\Z(2))\cong H^3(\kappa(v),\Q/\Z(2))=0
\end{equation}
as $\kappa(v)$ has cohomological dimension 2.
Therefore there are evaluation maps
\begin{equation}\label{eval}
Y(K_v)\times H^3(Y, \Q/\Z(2))\to H^3(K_v,\Q/\Z(2))\cong\Q/\Z,
\end{equation}
and moreover the similarly defined pairing
$$
{\mathcal Y}(\calo_v)\times H^3({\mathcal Y}, \Q/\Z(2))\to H^3(\calo_v,\Q/\Z(2))
$$
is trivial as its target is 0.
Hence the sum of the maps (\ref{eval}) defines a pairing
\begin{equation}\label{recpairing}
Y(\A_K)\times H^3(Y,\Q/\Z(2))\to\Q/\Z.
\end{equation}
The sequence
$$
H^3(K,\Q/\Z(2))\to\bigoplus_{v\in X^{(1)}} H^3(K_v,\Q/\Z(2))\stackrel\Sigma\longrightarrow \Q / \Z
$$
is a complex by virtue of the generalized Weil reciprocity law (\cite{cogal}, Chapter II, Annexe, (3.3) and (2.2)). Therefore the pairing (\ref{recpairing}) annihilates the diagonal image of $Y(K)$ in $Y(\A_K)$ on the left, and it also annihilates
the image of the natural map $H^3(K,\Q/\Z(2))\to H^3(Y,\Q/\Z(2))$
on the right. Set

\begin{align*}
H^3_{\rm lc}(Y,\Q/\Z(2)) &:= \ker(H^3(Y,\Q/\Z(2))/\im(H^3(K,\Q/\Z(2)))\to \\ &\to\prod_vH^3(Y\times_K{K_v},\Q/\Z(2))/\im(H^3(K_v,\Q/\Z(2)))).
\end{align*}
Then (\ref{recpairing}) induces a pairing
$$
[\,\,,\,\,]:\,Y(\A_K)\times H^3_{\rm lc}(Y,\Q/\Z(2))\to\Q/\Z
$$
annihilating the diagonal image of $Y(K)$ on the left. Since the elements of $H^3_{\rm lc}(Y,\Q/\Z(2))$ correspond to locally constant elements in $H^3(Y,\Q/\Z(2))$, for fixed $\alpha\in  H^3_{\rm lc}(Y,\Q/\Z(2))$ the value of $[(P_v),\alpha]$ does not depend on the choice of the adelic point $(P_v)\in Y(\A_K)$. Denote this common value by $\rho_Y(\alpha)$. We obtain a homomorphism
$$
\rho_Y:\,  H^3_{\rm lc}(Y,\Q/\Z(2))\to \Q/\Z.
$$

\begin{theo} \label{mainph}
Assume that $Y$ is a torsor under a $K$-torus $T$ such that $Y(\A_K)\neq\emptyset$ and $\rho_Y$ is the zero map. Then $Y(K)\neq\emptyset$.
\end{theo}

The proof follows a by now classical pattern (see e.g. \cite{skobook}, Theorems 6.2.1, 6.2.3), with additional complications caused by the fact that we are working over a field of cohomological dimension 3. We first construct a map
$$
\tau:\, \Sha^2(K, T')\to H^3_{\rm lc}(Y,\Q/\Z(2))
$$
whose role will be to relate the map $\rho_Y$ to the duality pairing of Theorem \ref{dualsha}. The map $\tau$ will be obtained from a map
\begin{equation}\label{H2map}
H^2(K, T')\to H^3(Y,\Q/\Z(2))/\im(H^3(K,\Q/\Z(2)))
\end{equation}
by passing to locally trivial elements. In order to define (\ref{H2map}), we consider the Hochschild--Serre spectral sequence
$$
H^p(K, H^q(\overline Y, \Q/\Z(2))\Rightarrow H^{p+q}(Y, \Q/\Z(2))
$$
where $\ov Y$ stands for $Y\times_K\ov K$ as usual. Since $K$ has cohomological dimension 3, outgoing differentials from the $E_2$-term $H^2(K, H^1(\overline Y, \Q/\Z(2))$ are trivial.
Therefore we have a map
\begin{equation}\label{preH2map}
H^2(K,H^1(\overline Y, \Q/\Z(2) )\to H^3(Y,\Q/\Z(2))/\im(H^3(K,\Q/\Z(2))).
\end{equation}
To complete the construction of the map (\ref{H2map}), we first show:

\begin{lem}\label{lemiso} There is an isomorphism of Galois modules
$$
H^1(\overline Y, \Q/\Z(2))\cong T'({\overline K})_{\rm tors}
$$
where $T'({\overline K})_{\rm tors}$ is the torsion submodule of $T'({\overline K})$.
\end{lem}

\begin{dem}
Fix an integer $n>0$. First of all, as a consequence of (\cite{sansuc}, Lemma 6.7) we have an isomorphism of Galois modules
$$
H^1(\overline Y, \mu_n)\cong H^1(\overline T, \mu_n).
$$
Since $H^1(\overline T,\G)=0$, the Kummer sequence in \'etale cohomology induces the first isomorphism in the chain
$$
H^1(\overline T, \mu_n)\cong H^0(\overline T, \G)/nH^0(\overline T, \G)\cong \widehat T/n\widehat T.
$$
The second one holds because by Rosenlicht's lemma $H^0(\overline T, \G)$ is an extension of $\widehat T$ by the divisible group $\overline K^\times$. Using the duality between $\widehat T$ and $\widecheck T$ we obtain an isomorphism of Galois modules
\begin{equation}\label{H1Y}
H^1(\overline Y, \mu_n)\cong \Hom_{\Z/n\Z}(\widecheck T/n\widecheck T, \Z/n\Z).
\end{equation}
Twisting by $\mu_n$ then yields
$$
H^1(\overline Y, \mu_n^{\otimes 2})\cong \Hom_{\Z/n\Z}(\widecheck T/n\widecheck T, \mu_n)\cong {}_nT'({\overline K})
$$
by definition of $T'$, whence the lemma follows by passing to the limit.
\end{dem}

By the lemma, we have a chain of isomorphisms
$$
H^2(K,H^1(\overline Y, \Q/\Z(2) )\cong H^2(K,T'({\overline K})_{\rm tors})\cong H^2(K,T')
$$
as the Galois module $T'({\overline K})/T'({\overline K})_{\rm tors}$ is uniquely divisible and therefore cohomologically trivial. This completes the construction of the map (\ref{H2map}), and hence also that of  $\tau$. \smallskip

We can now state:

\begin{prop}\label{mainphprop} The equality
$$
\rho_Y(\tau(\alpha))=\langle [Y],\alpha\rangle
$$
holds up to a sign, where $[Y]$ is the class of $Y$ in $\Sha^1(K,T)\subset H^1(K,T)$ and $\langle\,\,,\,\,\rangle$ is the duality pairing $\Sha^1(K,T)\times\Sha^2(K,T')\to\Q/\Z$ constructed in Theorem  \ref{dualsha}.
\end{prop}

Note that {\em Proposition \ref{mainphprop} implies Theorem \ref{mainph}.\/} Indeed, since the duality pairing of Theorem  \ref{dualsha} is non-degenerate, the assumption $\rho_Y=0$ in Theorem \ref{mainph} implies $[Y]=0$ in view of  the proposition, i.e. $Y(K)\neq \emptyset$. In fact, the argument shows that already the vanishing of $\rho_Y$ on $\im(\tau)\subset H^3_{\rm lc}(Y,\Q/\Z(2))$ implies $Y(K)\neq \emptyset$.  \enddem\smallskip

The rest of this section will be devoted to the {\em proof of Proposition \ref{mainphprop}.\/} We begin by another description of the map $\rho_Y$. Consider the commutative
diagram with exact rows
{\small $$
\begin{CD}
 H^3(K) @>>> H^3(Y) @>>> H^3(Y)/\im(H^3(K)) @>>> 0 \\
 @VVV @VVV @VVV \\
 \bigoplus_{v\in X^{(1)}} H^3(K_v) @>>> \bigoplus_{v\in X^{(1)}}H^3(Y\times_KK_v) @>>>
\bigoplus_{v\in X^{(1)}} H^3(Y\times_KK_v)/H^3(K_v) @>>> 0
\end{CD}
$$}
\noindent where all cohomology groups have $\Q/\Z(2)$-coefficients. The bottom left map is injective since $Y(\A_K)\neq \emptyset$. Therefore the snake lemma yields a map

{\small $$
{\rm Ker}[H^3(Y)/\im(H^3(K))\to\bigoplus_v H^3(Y\times_KK_v)/H^3(K_v)]\to {\rm Coker}(H^3(K)\to \bigoplus_{v} H^3(K_v)).
$$}

\noindent The kernel on the left here is none but $H^3_{\rm lc}(Y,\Q/\Z(2))$. The cokernel maps to $\Q/\Z$ by the reciprocity law recalled above. We thus obtain a map $H^3_{\rm lc}(Y,\Q/\Z(2))\to\Q/\Z$ which equals $\rho_Y$ by a checking analogous to the proof of (\cite{dhsza2}, Lemma 3.1).

Next, denote by  $\ov\pi$ the structure map $\ov Y\to\ov K$. Since $H^0(\ov Y,\Q/\Z(2))\cong\Q/\Z(2)$, we have a distinguished triangle in the bounded derived category of Galois modules over $K$ given by
\begin{equation}\label{triangle1}
\Q/\Z(2)\to\tau_{\leq 1} {\bf R}\ov\pi_*\Q/\Z(2)\to H^1(\ov Y,\Q/\Z(2))[-1]\to \Q/\Z(2)[1]
\end{equation}
where $\tau$ denotes the sophisticated truncation. It gives rise to the top row in the commutative exact diagram
$$
\begin{CD}
H^3(K, \Q/\Z(2)) @>>> H^3(K, \tau_{\leq 1}{\bf R}\ov\pi_*\Q/\Z(2)) @>>> H^2(K, H^1(\ov Y, \Q/\Z(2))) \\
@V{\id}VV @VVV @VVV \\
 H^3(K, \Q/\Z(2)) @>>> H^3(Y, \Q/\Z(2)) @>>> H^3(Y, \Q/\Z(2))/\im(H^3(K, \Q/\Z(2)))
\end{CD}
$$
where the right vertical map is given by (\ref{preH2map}) and the middle one by the composite
$$
H^3(K, \tau_{\leq 1}{\bf R}\ov\pi_*\Q/\Z(2))\to H^3(K, {\bf R}\ov\pi_*\Q/\Z(2))\cong H^3(Y, \Q/\Z(2)).
$$
Actually, the top right map in the above diagram is surjective because $K$ has cohomological dimension 3, and therefore the same snake-lemma construction as above (applied to the top row of the diagram over $K$ and its completions) yields a map
$$
\sigma_Y:\, \Sha^2(K, H^1(\ov Y, \Q/\Z(2)))\to \Q/\Z
$$
compatible with $\rho_Y$ via the map (\ref{preH2map}).

Here is another construction of the map $\sigma_Y$. Take $U\subset X$ so that there is a smooth model $\pi:\, {\mathcal Y}\to U$ of $Y$ over $U$ as above. Triangle (\ref{triangle1}) extends to a distinguished triangle in the bounded derived category of \'etale sheaves on $U$ given by
\begin{equation}\label{triangle2}
\Q/\Z(2)\to\tau_{\leq 1} {\bf R}\pi_*\Q/\Z(2)\to{\bf R}^1\pi_*\Q/\Z(2)[-1]\to \Q/\Z(2)[1]
\end{equation}
whence a map
$$
{\bf R}^1\pi_*\Q/\Z(2)\to \Q/\Z(2)[2]
$$
inducing
\begin{equation}\label{umap}
H^2_c(U, {\bf R}^1\pi_*\Q/\Z(2))\to H^4_c(U,\Q/\Z(2))\cong \Q/\Z
\end{equation}
where the last isomorphism was established in the course of the proof of Lemma \ref{globaltrace}. Since in view of Lemma \ref{lemiso} each class in $H^1(\ov Y,\Q/\Z(2))$ is represented by a point of the dual torus, a similar localization sequence as in Corollary \ref{hensel} shows that each class $\alpha\in \Sha^2(K, H^1(\ov Y, \Q/\Z(2)))$ lifts to a class $\alpha_U$ in $H^2_c(U, {\bf R}^1\pi_*\Q/\Z(2))$.

\begin{lem}
The image of $\alpha_U$ by the map $(\ref{umap})$ equals $\sigma_Y(\alpha)$.
\end{lem}

\begin{dem} For $v\in X^{(1)}$ denote by $j_v$ the natural map $\spec K_v^h\to U$, and consider the commutative diagram
$$
\begin{CD}
\Q/\Z(2)[1] @>>> \tau_{\leq 1} {\bf R}\pi_*\Q/\Z(2)[1] @>>> {\bf R}^1\pi_*\Q/\Z(2) \\
@VVV @VVV @VVV \\
\bigoplus_{v\notin U}j_{v*}j_v^*\Q/\Z(2)[1] @>>> \bigoplus_{v\notin
U}j_{v*}j_v^*\tau_{\leq 1} {\bf R}\pi_*\Q/\Z(2)[1] @>>> \bigoplus_{v\notin U}j_{v*}j_v^*{\bf R}^1\pi_*\Q/\Z(2)
\end{CD}
$$
whose rows are parts of exact triangles. Denote by $C_l$ the cone of the left vertical map and by $C_r$ that of the right one. Using Proposition \ref{commute} (1) we may identify $\alpha_U$ with an element in $H^2(C_r[-1])$, and its image by $(\ref{umap})$ with an element  $\beta_U\in H^2(C_l)$. Passing to the direct limit over $U$ smaller and smaller, the limit of the elements $\alpha_U$ becomes $\alpha$ and that of the $\beta_U$ the image of $\alpha$ by the snake-lemma construction. To see that the latter is indeed $\sigma_Y(\alpha)$, it remains to check that we may replace henselizations by completions in its construction. But the localization sequence in \'etale cohomology induces isomorphisms $H^3(K_v,\Q/\Z(2))\cong H^2(\kappa(v),\Q/\Z(1))$ in view of the vanishing of the groups in (\ref{vanish}), and the same argument yields isomorphisms $H^3(K_v^h,\Q/\Z(2))\cong H^2(\kappa(v),\Q/\Z(1))$.
\end{dem}

Let us now reexamine the duality pairing $\Sha^1(T)\times \Sha^2(T')\to\Q/\Z$ of Theorem \ref{dualsha}. According to its construction it is induced by pairings
\begin{equation}\label{cuppairing}
H^1(U, \T) \times H^2_c(U, \T')\to H^5_c(U,\Z(2))
\end{equation}
fitting into a commutative diagram

\vbox{$$
H^1(U, \T)\qquad  \times\qquad H^2_c(U, \T')
\quad\to\quad H^5_c(U,\Z(2))
$$

\begin{equation}\label{cupdiag}
\partial_n\downarrow\qquad\qquad\qquad \quad\quad
\uparrow\iota_n\qquad\quad\qquad \qquad\uparrow\partial
\end{equation}

$$
H^2(U, {}_n \T)\quad \times\qquad  H^2_c(U, {}_n\T')\to \quad H^4_c(U,
\Z/n\Z(2))
$$}

\noindent for all $n>0$, meaning $\alpha_1\cup \iota_n(\alpha_2)=\partial(\partial_n(\alpha_1)\cup\alpha_2)$ for $\alpha_1\in H^1(U,\T)$, ${\alpha_2\in  H^2_c(U, {}_n\T')}$. Indeed, rewriting the diagram as

\vbox{$$
H^1(U, \widecheck\T\otimes\Z(1)[1])\qquad  \times\qquad H^2_c(U, \widehat\T\otimes\Z(1)[1])
\quad\to\quad H^4_c(U,\Z(2)[1])
$$

$$
\partial_n\downarrow\qquad\qquad\qquad\qquad \qquad\qquad
\uparrow\iota_n\qquad\qquad\qquad \qquad\qquad\uparrow\partial
$$

$$
H^1(U, \widecheck\T\otimes\Z/n\Z(1)[1])\quad \times\qquad  H^2_c(U, \widehat\T\otimes\Z/n\Z(1))\to \quad H^4_c(U,
\Z/n\Z(2))
$$}

\noindent we see that the required compatibility follows from functoriality of the cup product pairing.
Since $H^2_c(U, \T')$ is torsion, each of its elements is of the form $\iota_n(\alpha_2)$ for suitable $n$. Thus the lower pairing in diagram (\ref{cupdiag}) determines the upper one.

On the other hand, the spectral sequence
$$
H^p(U, Ext^q(\widehat \T, \G))\Rightarrow \Ext^{p+q}_U(\widehat \T, \G)
$$
gives an edge map
\begin{equation}\label{map1}
H^1(U, \T)=H^1(U, Hom(\widehat \T, \G))\to \Ext^1_U( \widehat\T, \G)
\end{equation}
which may be composed with the boundary map
\begin{equation}\label{map2}
\Ext^1_U (\widehat\T, \G)\to \Ext^2_U (\widehat\T, \Z/n\Z(1)).
\end{equation}
Finally, tensoring morphisms in the derived category by $\Z/n\Z(1)$ yields a map
\begin{equation}\label{map3}
\Ext^2_U (\widehat\T, \Z/n\Z(1))\to \Ext^2_U (\widehat\T\otimes\Z/n\Z(1), \Z/n\Z(2)).
\end{equation}
But $\widehat\T\otimes\Z/n\Z(1)\cong {}_n\T'$, so there is a cup-product pairing
\begin{equation}\label{extpairing}
\Ext^2_U (\widehat\T\otimes\Z/n\Z(1), \Z/n\Z(2))\times H^2_c(U, {}_n\T')\to  H^4_c(U,
\Z/n\Z(2)).
\end{equation}
Pairings (\ref{cuppairing}) and (\ref{extpairing}) together with the composite
\begin{equation}\label{comp}
H^1(U, \T)\to \Ext^2_U (\widehat\T\otimes\Z/n\Z(1), \Z/n\Z(2))
\end{equation}
of the maps (\ref{map1}), (\ref{map2}) and (\ref{map3}) give rise to the diagram in the following lemma.
\begin{lem}
The diagram

\vbox{$$
{}\qquad H^1(U, \T)\qquad\quad  \times\qquad H^2_c(U, \T')
\quad\to\quad H^5_c(U,\Z(2))
$$

$$
{}\qquad\downarrow\qquad\qquad\quad \qquad\qquad
\uparrow\iota_n\qquad\quad\qquad \qquad\uparrow\partial
$$

$$
\Ext^2_U (\widehat\T\otimes\Z/n\Z(1), \Z/n\Z(2)) \times  H^2_c(U, {}_n\T')\to  H^4_c(U,
\Z/n\Z(2))
$$}

\noindent commutes.
\end{lem}

By the same argument as above, the lower pairing in the diagram of the lemma determines the upper one.\smallskip

\begin{dem} Recall that $\T\cong {Hom}(\widehat\T,\Z(1)[1])$ and $\T'\cong \widehat\T\otimes\Z(1)[1]$.  Tensoring by $\Z(1)[1]$ induces a map
${Hom}(\widehat\T,\Z(1)[1])\to {Hom}(\widehat\T\otimes\Z(1)[1],\Z(2)[2])$ making the diagram

\vbox{$$
{}\qquad H^1(U, \T) \qquad\qquad \times\qquad\qquad H^2_c(U, \T')\qquad\to\qquad
H^5_c(U,\Z(2))
$$

\begin{equation}\label{diag3}
{}\qquad\downarrow\qquad\qquad\quad \qquad\qquad
\qquad\qquad\downarrow\,\cong\qquad\quad\qquad\qquad\qquad \downarrow\id
\end{equation}

$$
H^1(U,{Hom}(\widehat\T\otimes\Z(1)[1],\Z(2)[2])) \times H^2_c(U, \widehat\T\otimes\Z(1)[1])\to H^5_c(U,
\Z(2))
$$}

\noindent commute. The edge map
\begin{equation}\label{edge}
H^1(U,{Hom}(\widehat\T\otimes\Z(1)[1],\Z(2)[2]))\to \Ext^1_U(\widehat\T\otimes\Z(1)[1],\Z(2)[2])
\end{equation}
in the Ext spectral sequence makes the lower pairing in this diagram compatible with the pairing
\begin{equation}\label{extpairing2}
\Ext^1_U(\widehat\T\otimes\Z(1)[1],\Z(2)[2]) \times H^2_c(U, \widehat\T\otimes\Z(1)[1])\to H^5_c(U,
\Z(2)).
\end{equation}
On the other hand, tensoring with $\Z/n\Z$ induces a map
\begin{equation}\label{tensn}
\Ext^1_U(\widehat\T\otimes\Z(1)[1],\Z(2)[2])\to \Ext^1_U(\widehat\T\otimes\Z/n\Z(1)[1],\Z/n\Z(2)[2])
\end{equation}
where the last group is none but $\Ext^2_U(\widehat\T\otimes\Z/n\Z(1),\Z/n\Z(2))$. We have a commutative diagram
$$
\begin{CD}
 H^1(U, \T) @>>> H^1(U,{Hom}(\widehat\T\otimes\Z(1)[1],\Z(2)[2])) \\
@VVV @VVV \\
\Ext^2_U(\widehat\T\otimes\Z/n\Z(1),\Z/n\Z(2)) @<<< \Ext^1_U(\widehat\T\otimes\Z(1)[1],\Z(2)[2])
\end{CD}
$$
whose arrows come clockwise from diagram (\ref{diag3}) and (\ref{edge}), (\ref{tensn}), (\ref{comp}), respectively.
Therefore the lemma reduces to the compatibility of the pairing (\ref{extpairing2}) with the lower pairing in the diagram of the lemma via the map (\ref{tensn}), which  is checked in the same way as the compatibility expressed by diagram (\ref{cupdiag}).
\end{dem}

We now relate the lower pairing in the lemma to the map (\ref{umap}). Each class in $\Ext^2_U (\widehat\T\otimes\Z/n\Z(1), \Z/n\Z(2))$ can be represented by a morphism $\widehat\T\otimes\Z/n\Z(1)\to \Z/n\Z(2)[2]$ in the bounded derived category of \'etale sheaves on $U$. But generalizing Lemma \ref{lemiso}, we have isomorphisms of \'etale sheaves
$$
\R^1\pi_*\Z/n\Z(2)\cong\widehat\T\otimes\Z/n\Z(1)\cong {}_n\T'
$$
(this uses the relative version of Rosenlicht's lemma; see \cite{desc}, Proposition 1.4.2). So passing to the direct limit over all $n$, we see that the pairing (\ref{extpairing}) induces a pairing
$$
\Ext^2_U (\R^1\pi_*\Q/\Z(2), \Q/\Z(2))\times H^2_c(U, \R^1\pi_*\Q/\Z(2))\to  H^4_c(U,
\Q/\Z(2)).
$$
A fixed morphism $\R^1\pi_*\Q/\Z(2)\to \Q/\Z(2)[2]$ in the derived category therefore induces a map
$$
H^2_c(U, \R^1\pi_*\Q/\Z(2))\to  H^4_c(U,
\Q/\Z(2))\cong\Q/\Z.
$$
This is exactly the same kind of map as (\ref{umap}) which is induced by the morphism $\R^1\pi_*\Q/\Z(2)\to \Q/\Z(2)[2]$ associated with $\tau_{\leq 1}\R\pi_*\Q/\Z(2)$. By passing to the direct limit over all $U$ we obtain the map
\begin{equation}\label{thismap}
H^1(\ov Y,\Q/\Z(2))\to \Q/\Z(2)[2]
\end{equation}
in the bounded derived category of Galois modules over $K$  that is associated with $\tau_{\leq 1}\R\ov\pi_*\Q/\Z(2)$ via the triangle (\ref{triangle1}). So to finish the proof of Proposition \ref{mainphprop}, and thereby that of Theorem \ref{mainph}, it suffices to check:

\begin{lem} The class of the map $(\ref{thismap})$ in $\Ext^2_K(\R^1\overline\pi_*\Q/\Z(2), \Q/\Z(2))$ is, up to a sign, the image of that of $[Y]\in H^1(K,T)$ via the composite map
$$
H^1(K, T)\to \Ext^2_K (\widehat T\otimes\Z/n\Z(1), \Z/n\Z(2))\to \Ext^2_K(H^1(\ov Y,\Q/\Z(2)), \Q/\Z(2)),
$$
where the first map is obtained from $(\ref{comp})$ by passing to the direct limit over $U$ and the second one comes from Lemma \ref{lemiso}.
\end{lem}

\begin{dem} It is enough to check that the image of $[Y]$ by the composite map
$$
H^1(K,T)\to \Ext^1_K (\widehat T, \G)\to \Ext^2_K (\widehat T, \Z/n\Z(1))\to \Ext^2_K(H^1(\ov Y,\Z/n\Z(1)), \Z/n\Z(1))
$$
corresponds up to a sign to the map $H^1(\ov Y,\Z/n\Z(1))\to \Z/n\Z(1)[2]$ coming from the triangle
$$
\Z/n\Z(1)\to\tau_{\leq 1}\R\ov\pi_*\Z/n\Z(1)\to H^1(\ov Y,\Z/n\Z(1))[-1]\to \Z/n\Z(1)[1]. $$
But it is well known (see e.g. \cite{skobook}, Lemma 2.4.3) that the class $[Y]$ corresponds, up to a sign, to the class in $\Ext^1_K (\widehat T, \G)$ of the extension
$$
1\to\G\to \pi_*\G\to\widehat T\to 1
$$
coming from Rosenlicht's lemma. On the other hand, it follows from (\cite{bvh}, Lemma 2.3, Corollary 2.5 and Lemma 4.2) that the same class is also represented by $\tau_{\leq 1}\R\ov\pi_*\G$. Now the lemma follows by composing with the map $\G\to\Z/n\Z(1)[1]$.

\end{dem}

We close this section by the following statement related to Theorem \ref{mainph} pointed out to us by J-L. Colliot-Th\'el\`ene.

\begin{cor}
Assume that $Y$ is a torsor under a torus $T$ that is stably rational over $K$.  If $Y(\A_K)\neq\emptyset$, then $Y(K)\neq\emptyset$.
\end{cor}

\begin{dem}
The point is that for a stably rational $K$-torus we have $\Sha^2(T')=0$, which implies $\Sha^1(T)=0$ by Theorem \ref{dualsha}. Alternatively, it implies that $\rho_Y$ vanishes on the image of $\Sha^2(T')$ in $H^3_{\rm lc}(Y,\Q/\Z(2))$, so the refined form of Theorem \ref{mainph} mentioned after Proposition \ref{mainphprop} applies.

To show the vanishing of $\Sha^2(T')$, one uses a theorem of Voskresenski{\u \i} \cite{vosk1}: for $T$ stably $K$-rational there is an exact sequence
$$
0\to \widehat T\to\widehat Q\to\widehat F\to 0
$$
where both $\widehat Q$ and $\widehat F$ are permutation modules (see also \cite{requiv}, Proposition 6). Tensoring by $\G$ we obtain a resolution of the dual torus $T'$ by the quasi-trivial tori $Q'$ and $F'$ (note that the dual of a quasi-trivial torus is again quasi-trivial). The resolution induces a map $\Sha^2(T')\to \Sha^2(Q')$ which is injective since $H^1(K,F')=0$ by Hilbert's Theorem 90. But $\Sha^2(Q')=0$  by Remark \ref{shaq}.
\end{dem}

\section{A generalization to reductive linear algebraic groups}\label{last}

In this section we prove a generalization of Theorem~\ref{mainph} to a class of noncommutative linear algebraic groups. The price we have to pay is a good reduction assumption. Therefore throughout the whole section {\it we assume that $K$ is the function field of a smooth projective curve $X$ defined over a p-adic field $k$  that extends over the ring of integers $\calo_k$ to
a proper and smooth relative curve $\X \to \spec \calo_k$ with geometrically
integral special fibre.}

Before proceeding further, we recall some standard notation and terminology concerning linear algebraic groups. Given a reductive group $G$, its derived subgroup is closed, connected and semisimple; we denote it by $G^{\rm ss}$. The semisimple group $G^{\rm ss}$ has a finite cover by a simply connected group $G^{\rm sc}$ which is unique up to isomorphism. The group $G$ is said to be  quasi-split if it has a Borel subgroup defined over the base field. Finally, $G$ is absolutely almost simple if over the algebraic closure it becomes an extension of a simple group by a finite normal (hence central) subgroup. \smallskip

Our main result concerns the reciprocity obstruction ${\rho_E: H^3_{\rm lc}(E, \Q/\Z(2))\to\Q/\Z}$  introduced in Section \ref{five}.

\begin{theo}\label{phlag}
Let $E$ be a torsor under a reductive linear algebraic group $G$ over $K$. Assume that $K$ satisfies the good reduction assumption made above and that $G^{\rm sc}$ is quasi-split and has no $E_8$ factor. If $E(\A_K)\neq\emptyset$ and  $\rho_E: H^3_{\rm lc}(E, \Q/\Z(2))\to\Q/\Z$ is the zero map, then $E(K)\neq\emptyset$.
\end{theo}

There are two main inputs for the proof. One is Theorem \ref{mainph}, the other is the following variant of Theorem 5.4 of \cite{ctps}.

\begin{prop} \label{sssc}
Let $G$ be a quasi-split, simply connected semisimple group over $K$ as above that has no $E_8$ factor. The kernel of the map
$$H^1(K,G) \to \prod_{v \in X^{(1)}} H^1(K_v,G)$$
is trivial.
\end{prop}

\begin{dem}
Assume first that $G$ is absolutely almost simple. Consider the commutative diagram
$$
\begin{CD}
H^1(K,G) @>>> \prod_{v \in X^{(1)}} H^1(K_v,G) \\
@VVV @VVV \\
H^3(K,\Q/\Z(2)) @>>> \prod_{v \in X^{(1)}}
H^3(K_v,\Q/\Z(2))
\end{CD}
$$
whose vertical maps are induced by the Rost invariants $H^1(K, G)\to H^3(K, \Q/\Z(2))$  and $H^1(K_v, G)\to H^3(K_v, \Q/\Z(2))$. As proven in (\cite{ctps}, Theorem 5.3), these maps have trivial kernel for $G$ as in the statement because $K$ and the $K_v$ have cohomological dimension 3. It therefore suffices to prove that the lower horizontal map has trivial kernel. According to Kato (\cite{katocrelle}, Prop. 5.2),
this kernel is isomorphic to the kernel of the
map
$$H^2(\kappa,\Q/\Z(1)) \to \prod_{\lambda \in X_0^{(1)} }
H^2(\kappa_{\lambda},\Q/\Z(1)),$$
where $\kappa$ is the function field of the special fibre $X_0$ of $X$ and the $\kappa_\lambda$ are its completions at the codimension 1 points $\lambda$. But by assumption $X_0$ is smooth, and therefore the above map is injective by the classical Brauer-Hasse-Noether theorem over function fields.

In the general case $G$
is a direct product of finitely many $G_i$ of the form
${G_i=R_{K_i/K}(H_i)}$, with $H_i$ simply connected and absolutely
almost simple
over $K_i$ (see \cite{kneser}, Hilfssatz 7.4, 7.5). Observe that $H_i$ is also quasi-split and without $E_8$
factor since it is a direct factor of $G_i \times_K K_i$. We have proven above the lemma for the $H_i$, so it also holds for the $G_i$ by the noncommutative Shapiro lemma (\cite{inv}, Lemma 29.6). Finally, we obtain the result for $G$ by taking the direct product.
\end{dem}

\begin{remas}\rm ${}$

\noindent (1) As of today, no example seems to be known where the kernel in the proposition is nontrivial for a simply connected $G$.

\noindent (2) In (\cite{ctps}, Theorem 5.5) triviality of the kernel of the Rost invariant is also proven for groups of $E_8$ type, provided that $G$ is defined over the ring of integers $\calo_k$ of $k$ and that the residue characteristic is at least 5. So the proposition (and with it the theorem) also hold in this case.

Also, Mathieu Florence and Philippe Gille point out to us that if $G$ is of type $E_8$, it follows from results of Chernousov (\cite{ch1}, \cite{ch2}) and Semenov \cite{semenov} that over a field of cohomological dimension  $\leq 4$ a $G$-torsor representing a class in the kernel of the Rost invariant has a zero-cycle of degree 1. So for groups of $E_8$ type Theorem \ref{phlag} will hold with this weaker conclusion (i.e. that $E$ as in the theorem has a zero-cycle of degree 1).

\noindent (3) The triviality of the kernel of the Rost invariant  used in the above proof is also known for some non-quasisplit groups over fields of cohomological dimension $\leq 3$ (\cite{preeti}, see also \cite{hu}). However, quasi-splitness will also play a crucial role in the proof of Proposition \ref{lemgt} below.
\end{remas}

To reduce the theorem to the cases of simply connected groups and tori, we need the construction of $z$-extensions first used by Deligne and Langlands (see \cite{mishih}, Prop. 3.1). It gives us a central extension of $K$-groups
\begin{equation}\label{zext}
1\to Q\to\widetilde G\to G\to 1
\end{equation}
where $Q$ is a quasi-trivial torus, and the derived subgroup of $\widetilde G$ is the simply connected cover $G^{\rm sc}$ of the derived subgroup $G^{\rm ss}$ of $G$.

\begin{lem}\label{isosha}
The above exact sequence induces an isomorphism of pointed sets $\Sha^1(\widetilde G)\stackrel\sim\to\Sha^1(G)$.
\end{lem}

\begin{dem}
The long exact cohomology sequences associated with the above central extension give rise to a commutative diagram of pointed sets with exact rows
{\small $$
\begin{CD}
H^1(K, Q)  @>>> H^1(K, \widetilde G) @>>> H^1(K, G) @>>> H^2(K, Q) \\
@VVV  @VVV @VVV @VVV \\
\prod_v H^1(K_v, Q) @>>>\prod_v H^1(K_v, \widetilde G) @>>>\prod_v H^1(K_v, G) @>>> \prod_v H^2(K_v, Q).
\end{CD}
$$}

\noindent Notice that $H^1(K, Q)=H^1(K_v, Q)=0$ for all $v$ as $H^1$ of a quasi-trivial torus is trivial over any field (a consequence of Hilbert's Theorem 90). Since $Q$ is central in $\widetilde G$, this shows that the middle horizontal maps are actually {\em injective} (see \cite{cogal}, Chapter I, Proposition 42).  On the other hand, Remark \ref{shaq}  implies that the right vertical map is injective as well. The proposition follows by a diagram chase.
\end{dem}

Since $\widetilde G$ is reductive, we have an exact sequence
$$
1\to G^{\rm sc}\to\widetilde G\to T\to 1
$$
with a torus $T$.

\begin{prop}\label{lemgt}
The map $\Sha^1(\widetilde G)\rightarrow\Sha^1(T)$ induced by the projection $\widetilde G\to T$  has trivial kernel.
\end{prop}

For the proof of the proposition we need two general statements about reductive groups for which we could not find an adequate reference.

\begin{lem}\label{maxtor}
Let $\widetilde G$ be a reductive group over a field with derived subgroup $G^{\rm ss}$. Every maximal torus of $\widetilde G$ is an extension of the torus  $T:=\widetilde G/G^{\rm ss}$ by a maximal torus of $G^{\rm ss}$ and conversely, such extensions yield maximal tori in $\widetilde G$.
\end{lem}

We thank Michel Brion for his help with the proof below.\medskip

\begin{dem}
Let $\widetilde T$ be a maximal torus in $\widetilde G$. We show that the map $\widetilde T\to T$ is surjective and that $\widetilde T\cap G^{\rm ss}$ is a maximal torus in $G^{\rm ss}$. Since every maximal torus of $G^{\rm ss}$  is contained in a maximal torus of $\widetilde G$, the lemma will follow from these facts.

To do so, we may assume
 the base field is algebraically closed. Since $\widetilde G$ is reductive, its radical $R(\widetilde G)$ has finite intersection with $G^{\rm ss}$ (\cite{springer}, Corollary 7.3.1 (ii)), and we have $\widetilde G= G^{\rm ss}\cdot R(\widetilde G)$ (\cite{springer}, Corollary 8.1.6). Therefore $R(\widetilde G)$ surjects onto $T$ via the projection $\widetilde G\to T$ with finite kernel. Furthermore, $R(\widetilde G)$ is central in $\widetilde G$ (\cite{springer}, Corollary 7.3.1 (i)) and is therefore contained in every maximal torus. As a consequence, every maximal torus of $\widetilde G$ surjects onto $T$ via the projection $\widetilde G\to T$, and also onto a maximal torus of $\widetilde G/R(\widetilde G)$ via the projection  $\widetilde G\to\widetilde G/R(\widetilde G)$.

 Consider now the identity component $S$ of $\widetilde T\cap G^{\rm ss}$. By the results of the previous paragraph, it is a torus of dimension $\dim \widetilde T-\dim T=\dim \widetilde T-\dim R(\widetilde G)$, which is the dimension of a maximal torus in $G^{\rm ss}$ (because the map $G^{\rm ss}\to\widetilde G/R(\widetilde G)$ is surjective with finite kernel). Therefore $S$ is a maximal torus in $G^{\rm ss}$.
 Now $\widetilde T\cap G^{\rm ss}$ certainly centralizes $S$ in $G^{\rm ss}$, but a maximal torus in a reductive group is its own centralizer (\cite{springer}, Corollary 7.6.4 (ii)). Hence indeed $S=\widetilde  T\cap G^{\rm ss}$.
\end{dem}

The following statement was pointed out to us by Jean-Louis Colliot-Th\'el\`ene.

\begin{lem}\label{quasisplit}
   A quasi-split and simply connected semi-simple algebraic group $G^{\rm sc}$ over an infinite field $F$ contains a quasi-trivial maximal torus defined over $F$.
\end{lem}

\begin{dem}
By quasi-splitness there is a Borel subgroup $B$ in $G^{\rm sc}$ defined over $F$. Since $F$ is infinite, $B$ contains a maximal torus $T_B$ defined over $F$ (\cite{springer}, Theorem 13.3.6). Over an algebraic closure $\overline F$ of $F$ the fibration of $F$-varieties $G\to G/B$ gives rise to an exact sequence of $\gal(\overline F|F)$-modules
$$
\widehat G^{\rm sc}\to \widehat B\to \pic (\overline {G^{\rm sc}/B})\to\pic \overline{G^{\rm sc}}
$$
by (\cite{sansuc}, Proposition 6.10 and Lemme 6.5 (iii)). Since $B$ is connected and solvable, we have an isomorphism of character groups $\widehat B\cong \widehat T_B$. Moreover, $\widehat G^{\rm sc}=0$ as $G^{\rm sc}$ is semisimple (\cite{springer}, Theorem 8.1.5 (ii)), and $\pic \overline {G^{\rm sc}}=0$ as $G^{\rm sc}$ is also simply connected (\cite{sansuc}, Lemme 6.5 (iv)). Now the lemma follows from the fact that $\pic (\overline {G^{\rm sc}/B})$ is a permutation module (see \cite{ctgp}, proof of Lemma 5.6 which is based on a result of Harder \cite{harder}).
\end{dem}

\noindent{\em Proof of Proposition \ref{lemgt}.}
 Each element of $\ker(\Sha^1(\widetilde G)\to\Sha^1(T))$ comes from an element of $H^1(K, G^{\rm sc})$. We show that this element actually lies in $\Sha^1(G^{\rm sc})$, i.e. it is locally trivial. This will conclude the proof in view of Proposition \ref{sssc}.

To show the required local triviality, we verify that for all completions $K_v$ the map $H^1(K_v, G^{\rm sc}) \to H^1(K_v,\widetilde G)$ has trivial kernel. Using the two previous lemmas and the assumption that $G^{\rm sc}$ is quasi-split we find an exact sequence
$$
0\to T^{\rm sc}\to\widetilde T\to T\to 0
$$
where $T^{\rm sc}$ is a quasi-trivial maximal torus in $G^{\rm sc}$ and $\widetilde T$ is a maximal torus in $\widetilde G$.
We have a commutative diagram with exact rows
$$
\begin{CD}
H^0(K_v, T) @>>> H^1(K_v, T^{\rm sc}) @>>> H^1(K_v,\widetilde T)  \\
@VV{\id}V @VVV @VVV  \\
H^0(K_v, T) @>>> H^1(K_v, G^{\rm sc}) @>>> H^1(K_v, \widetilde G).
\end{CD}
$$
Notice that  $H^1(K_v, T^{\rm sc})=1$ since $T^{\rm sc}$ is a quasi-trivial torus. The triviality of $\ker(H^1(K_v, G^{\rm sc}) \to H^1(K_v, \widetilde G))$ then follows by a diagram chase.
\enddem

\noindent{\em Proof of Theorem \ref{phlag}.} The proof is particularly simple when ${G=\widetilde  G}$, i.e. when $G^{\rm ss}$ is a simply connected group satisfying the assumptions of the theorem, so we first treat this case. Assume $E$ is a torsor under $G$ such that $E(\A_K)\neq\emptyset$. Denote by $\rho_E:\, H^3_{\rm lc}(E, \Q/\Z(2))\to\Q/\Z$ the map given by evaluating classes in $H^3_{\rm lc}(E, \Q/\Z(2))$ at adelic points of $E$. Denoting the torus $G/G^{\rm ss}$ by $T$ as above, the pushforward of $E$ by the map $G\to T$ gives a torsor $Y$ under $T$ with associated evaluation map ${\rho_Y: H^3_{\rm lc}(Y, \Q/\Z(2))\to\Q/\Z}$. The maps $\rho_Y$ and $\rho_E$ are compatible via the map
$$
H^3_{\rm lc}(Y, \Q/\Z(2))\to H^3_{\rm lc}(E, \Q/\Z(2))
$$
induced by $E\to Y$. Therefore if $\rho_E$ is trivial, so is $\rho_Y$. But by Theorem \ref{mainph} the torsor $Y$ is then trivial. By Proposition \ref{lemgt}, so is $E$.

To treat the general case, we have to delve into the constructions proving Theorem \ref{mainph}. We also need a definition of the algebraic fundamental group of $G$; for ease of reference we adopt that of \cite{ctresflasq}, \S 6 (which coincides with that of \cite{borams}, \S 1). It gives us a Galois module $\pi_1^{\rm alg}(G)$ which is of finite type as a $\Z$-module and coincides with the cocharacter group of $G$ in case $G$ is a torus. We do not need the precise definition, only the following properties: $\pi_1^{\rm alg}$ is an exact functor from the category of reductive groups over $K$ to that of Galois modules (\cite{ctresflasq}, Proposition 6.8), and for all $n>0$ there is an isomorphism of Galois modules
$$
H^1(\overline G,\mu_n)\cong \Hom_{\Z/n\Z}(\pi_1^{\rm alg}(G)/n, \Z/n\Z)
$$
(\cite{ctresflasq}, Proposition 6.7).

Comparing with the isomorphism (\ref{H1Y}) of Section \ref{five}, we see that by proceeding as in the construction of the map (\ref{H2map}) we obtain a map
$$
H^2(K, G^*)\to H^3(E, \Q/\Z(2))/\im(H^3(K,\Q/\Z(2))
$$
where $G^*$ is the group of multiplicative type whose character module is $\pi_1^{\rm alg}(G)$, and $E$ is a $K$-torsor under $G$. Passing to locally trivial elements yields a map
$$
\Sha^2(G^*)\to H^3_{\rm lc}(E, \Q/\Z(2)).
$$
By exactness of the functor $\pi_1^{\rm alg}$, exact sequence (\ref{zext}) gives rise to an exact sequence of Galois modules
$$
0\to \widecheck Q\to \pi_1^{\rm alg}(\widetilde G)\to \pi_1^{\rm alg}(G)\to 0,
$$
whence a dual exact sequence of groups of multiplicative type
\begin{equation}\label{qex}
1\to G^*\to \widetilde G^*\to Q'\to 1.
\end{equation}
Since $Q'$ is a quasi-trivial torus, we have $H^1(K, Q')=0$ and also $\Sha^2(Q')=0$ by Remark \ref{shaq}. Therefore the long exact cohomology sequence of (\ref{qex}) induces an isomorphism  $\Sha^2(G^*)\stackrel\sim\to \Sha^2(\widetilde G^*)$.

Now let again $E$ be a locally trivial $K$-torsor under $G$, and let $\widetilde E$ be a locally trivial $K$-torsor under $\widetilde G$ whose pushforward is $E$; such an $\widetilde E$ exists by Lemma \ref{isosha}. As in the first part of this proof, let $Y$ be the pushforward of $\widetilde E$ by the map $\widetilde G\to T$, where $T=\widetilde G/\widetilde G^{\rm ss}$. We have a commutative diagram
$$
\begin{CD}
\Sha^2(T') @>>> \Sha^2(\widetilde G^*) @<{\cong}<< \Sha^2(G^*) \\
@VVV @VVV @VVV \\
H^3_{\rm lc}(Y, \Q/\Z(2)) @>>> H^3_{\rm lc}(\widetilde E, \Q/\Z(2)) @<<< H^3_{\rm lc}(E, \Q/\Z(2)).
\end{CD}
$$
 A diagram chase shows that if the evaluation map $\rho_E:H^3_{\rm lc}(E, \Q/\Z(2))\to\Q/\Z$ is trivial, then so is the restriction of the evaluation map $\rho_Y:H^3_{\rm lc}(Y, \Q/\Z(2))\to\Q/\Z$ to the image of $\Sha^2(T')$ in $H^3_{\rm lc}(Y, \Q/\Z(2))$. By the refined form  of Theorem \ref{mainph} mentioned after Proposition \ref{mainphprop}, the triviality of this restriction already implies the triviality of $Y$, whence the triviality of $E$ by Lemma \ref{isosha} and Proposition \ref{lemgt}.\enddem

\end{document}